%
%
%
%
\ifx\ENSMATH\undefined\let\ENSMATH=Y\else  \fi
\input amstex
\UseAMSsymbols
%
%
\topskip=9.0mm  \advance \topskip by -14pt 
\hsize=11.3cm    
\vsize=17.6cm    
\hoffset=7.95truecm
\advance\hoffset by -5.65cm
\voffset=11.8truecm
\advance\voffset by -8.8cm
\parindent=.5cm     
\def\footnoterule{\kern -3pt\hrule width 1.4cm\kern 4pt}
\parskip=.5pt plus.75pt minus.25pt   
\hfuzz=2pt          
\mathsurround=.1em  
\pretolerance=400   
\tolerance=800      
\binoppenalty=2400  
\relpenalty=1200    
\clubpenalty=1000   
\widowpenalty=1000  
\frenchspacing      
\newskip\EMskip     
\EMskip=\normalbaselineskip
\def\EMquad{\hskip 1em plus .02em minus .2em}  
%
%
\abovedisplayskip=7pt plus 3pt minus 4pt
\belowdisplayskip=7pt plus 3pt minus 4pt
\abovedisplayshortskip=0pt plus 3pt
\belowdisplayshortskip=5pt plus 3pt minus 2pt
\predisplaypenalty=0
\def\EMpar{\belowdisplayskip=0pt\belowdisplayshortskip=0pt\par}
%
%
\catcode`\@=11
\catcode`\;=\active
\def;{\relax\ifhmode\ifdim\lastskip>\z@\unskip\fi\kern.1em
\fi\string;\ifhmode\ \ignorespaces\fi}
\catcode`\:=\active
\def:{\relax\ifhmode\ifdim\lastskip>\z@\unskip\fi\kern.15em
\fi\string:\ifhmode\ \ignorespaces\fi}

\let\colon=\:
\catcode`\!=\active
\def!{\relax\ifhmode\ifdim\lastskip>\z@\unskip\fi\kern.18em\fi\string!}

\catcode`\?=\active
\def?{\relax\ifhmode\ifdim\lastskip>\z@\unskip\fi\kern.18em\fi\string?}

%
%
\def\EMraise#1{\leavevmode\raise.78ex\hbox{#1}}
\def\EMtvn{\vrule height 2.9mm depth 0pt width 0pt}
\def\EMmark#1{\def\next{#1}\ifx\next\empty\else
\smash{\ninepol\EMraise{\sixpol #1}\kern.12em)}\fi\EMtvn}
\long\def\EMtext#1{\pluspetit\ignorespaces #1\par}
\long\def\footnote#1#2{\relax
\ifhmode\ifdim\lastskip>\z@\unskip\fi\/\kern.15em%
\bgroup\parindent=.8cm\plainfootnote{\bgroup\EMmark{#1}\egroup}%
{{\EMtext{#2}}\vskip -\baselineskip}\egroup\fi}
\newcount\notenumber  \notenumber=1
\long\def\note#1{\footnote{\the\notenumber}{#1}%
{\global\advance\notenumber by 1 }}
\footline={\hfil\tenpol\folio\hfil}
%
%
\let\maj=\uppercase

\let\emptyset=\varnothing

\def\bigpenalty{\interlinepenalty=\@M }
\def\smallitem#1{\par\vskip 2pt\noindent\hbox to .6cm{\ignorespaces 
#1\hfill}\hangindent=.6cm\hangafter=1 \ignorespaces }
\def\meditem#1{\par\vskip 2pt\noindent\hbox to .7cm{\hfill
\ignorespaces #1\enspace}%
\hangindent=.7cm\hangafter=1 \ignorespaces }
\def\bigitem#1{\par\vskip 2pt\noindent\hbox to .8cm{\hfill
\ignorespaces #1\enspace}%
\hangindent=.8cm\hangafter=1 \ignorespaces }
\def\cite#1{[#1]}
\let\INDENT=\indent
\def\CASE{\relax\ifhmode\ifdim\lastskip>\z@\unskip\fi\fi}
\newif\ifperiod  \periodtrue  
\newcount\LANGUE
\def\anglais{\LANGUE=0 }

\anglais
\def\BY{\ifcase\LANGUE by\or par\or von\or da\else\fi}
\def\AND{\ifcase\LANGUE and\or et\or und\or e\else\fi}
\def\ABSTRACT{\ifcase\LANGUE Abstract\or R\'esum\'e\or 
     Kurzfassung\or Riassunto\else\fi}
\def\DEFINITION{\ifcase\LANGUE Definition\or D\'efinition\or 
     Definition\or Definizione\else\fi}
\def\DEFINITIONS{\ifcase\LANGUE Definitions\or D\'efinitions\or 
     Definitionen\or Definizioni\else\fi}
\def\EXAMPLE{\ifcase\LANGUE Example\or Exemple\or 
     Beispiel\or Esempio\else\fi}
\def\EXAMPLES{\ifcase\LANGUE Examples\or Exemples\or 
     Beispiele\or Esempi\else\fi}
\def\PROOF{\ifcase\LANGUE Proof\or D\'emonstration\or 
     Beweis\or Dimostrazione\else\fi}
%
%
\newcount\FIRST
\newbox\EMboxone \newbox\EMboxtwo \newbox\EMboxthree
\newbox\EMboxfour \newbox\EMboxfive
\def\raggedcenter{\leftskip=0pt plus 1fill \rightskip=0pt plus 1fill
\parfillskip=0pt\pretolerance=\@M \hyphenpenalty=\@M }
\long\def\boitecentree#1{\bgroup\parindent=.8cm
\bgroup\global\FIRST=0
\long\def\footnote##1##2{\relax                               
\ifhmode\ifdim\lastskip>\z@\unskip\fi\/\kern.15em\fi%
\global\advance\FIRST by 1
\ifnum\FIRST=1 \global\setbox\EMboxone=\hbox{\EMmark{##1}}%
\gdef\EMtextone{\EMtext{##2}}\copy\EMboxone\fi%
\ifnum\FIRST=2 \global\setbox\EMboxtwo=\hbox{\EMmark{##1}}%
\gdef\EMtexttwo{\EMtext{##2}}\copy\EMboxtwo\fi%
\ifnum\FIRST=3 \global\setbox\EMboxthree=\hbox{\EMmark{##1}}%
\gdef\EMtextthree{\EMtext{##2}}\copy\EMboxthree\fi%
\ifnum\FIRST=4 \global\setbox\EMboxfour=\hbox{\EMmark{##1}}%
\gdef\EMtextfour{\EMtext{##2}}\copy\EMboxfour\fi%
\ifnum\FIRST=5 \global\setbox\EMboxfive=\hbox{\EMmark{##1}}%
\gdef\EMtextfive{\EMtext{##2}}\copy\EMboxfive\fi}
\vbox{\parindent=0pt\raggedcenter\let\par=\endgraf\let\\=\par
#1\par}\egroup%
\ifnum\FIRST>0 \vfootnote{\copy\EMboxone}{{\EMtextone}%
\vskip -\baselineskip}\fi
\ifnum\FIRST>1 \vfootnote{\copy\EMboxtwo}{{\EMtexttwo}%
\vskip -\baselineskip}\fi
\ifnum\FIRST>2 \vfootnote{\copy\EMboxthree}{{\EMtextthree}%
\vskip -\baselineskip}\fi
\ifnum\FIRST>3 \vfootnote{\copy\EMboxfour}{{\EMtextfour}%
\vskip -\baselineskip}\fi
\ifnum\FIRST>4 \vfootnote{\copy\EMboxfive}{{\EMtextfive}%
\vskip -\baselineskip}\fi
\egroup}
%
%
\newtoks\titrecrt  \titrecrt={\hfil}
\long\def\titre#1{\null\vskip 2.5cm plus 1.2cm minus .5cm
\boitecentree{\let\desserre=\desserreplus\normal\baselineskip=5mm #1}%
\ifx\TITREBREF\undefined\titrecrt={%
\def\\{\CASE\space\ignorespaces }\let\par=\\%
#1}\fi\vskip 1cm}
\def\titrebref#1{\let\TITREBREF=Y\titrecrt={#1}}
%
%
\newtoks\auteurcrt  \auteurcrt={\hfil}
\long\def\auteur#1{\vskip -.6cm
\def\sn##1{\bgroup\smc ##1\egroup}%
\boitecentree{\normal\BY\ {#1}}%
\ifx\AUTEURBREF\undefined\auteurcrt={%
\long\def\footnote##1##2{\relax}\def\sn##1{##1}%
\def\\{\CASE\space\ignorespaces }\let\par=\\%
\maj{#1}}\fi\vskip 1cm}
\def\auteurbref#1{\let\AUTEURBREF=Y\auteurcrt={\maj{#1}}}
%
%
\long\def\EMabstract#1#2{\goodbreak\bgroup\petit\INDENT
\def\next{#1}\ifx\next\empty\else{\smc {\ignorespaces #1\CASE}.\EMquad}\fi
\rm\ignorespaces #2\EMpar\egroup\vskip .25cm}
\long\def\abstract#1{\EMabstract{\ABSTRACT}{#1}}

%
%
\long\def\sec#1#2{\removelastskip\goodbreak\vskip 10mm plus 6mm minus 3mm
\normal\goodbreak\boitecentree{\let\it=\sl
\smc\def\next{#1}\ifx\next\empty\ignorespaces #2\else
\def\\##1/{\labchap{##1}.}%
\ignorespaces #1\EMquad\ignorespaces #2\fi}
\nobreak\vskip\baselineskip}
\long\def\subsec#1#2{\removelastskip\penalty-200 \vskip\EMskip
\penalty-200 \noindent{\let\it=\sl
\smc\def\next{#1}\noindent\ifx\next\empty\ignorespaces #2\else
\def\\##1/{\labsec{##1}}%
\setbox0=\hbox{\ignorespaces #1\EMquad}%
\dimen1=\wd0\box0\hangindent=\dimen1\hangafter=1 
\def\\{\CASE\break\ignorespaces}%
\ignorespaces #2\fi\par}\par\nobreak\vskip.35\baselineskip}
\long\def\parag#1{\removelastskip\penalty-100 \vskip\EMskip\INDENT
\leavevmode\null\def\next{#1}\ifx\next\empty\else{%
\def\\##1/{\labsec{##1}}%
\smc {\ignorespaces #1\CASE}\ifperiod.\else
\global\periodtrue\fi\kern.8em}\fi\ignorespaces }
%
%
%
\long\def\proclaim#1#2{\ifvmode\removelastskip\goodbreak\vskip\EMskip
\INDENT\fi\def\next{#1}\ifx\next\empty\else{%
\def\\##1/{\labtheo{##1}}%
\smc {\ignorespaces
#1\CASE}\ifperiod.\else\global\periodtrue\fi\EMquad
}\fi{\interlinepenalty=\@M \clubpenalty=5000
\widowpenalty=\@M \predisplaypenalty=\@M
\it\ignorespaces #2\EMpar}\vskip\EMskip}
%
%




\long\def\lemma#1#2{\proclaim{Lemma~\rm\ignorespaces #1}{#2}}




\long\def\theorem#1#2{\proclaim{Theorem~\rm\ignorespaces #1}{#2}}

%
%
\long\def\EMdemo#1#2{\removelastskip\vskip.5\baselineskip
\INDENT\def\next{#1}\ifx\next\empty\else{\it {\ignorespaces 
#1\CASE}\ifperiod.\else\global\periodtrue\fi\EMquad}\fi{\clubpenalty=5000
\rm\ignorespaces #2\EMpar}\penalty-150\vskip\EMskip}

\long\def\proof#1{\EMdemo{\PROOF}{#1}}

%
%
\newbox\sqwh    
\def\makesqwh{\def\EMtraiv{\vrule height 1.5ex depth .03ex width .03em}%
\def\EMtraih{\hrule height .03em depth 0pt}%
\setbox\sqwh=\vbox{\EMtraih\hbox{\EMtraiv\kern 1.53ex\EMtraiv}\EMtraih}%
\dimen1=1.5ex  \advance\dimen1 by .03em  \ht\sqwh=\dimen1   
\dimen2=.03ex  \advance\dimen2 by .03em  \dp\sqwh=\dimen2 \relax}
\makesqwh
\def\qed{\ifmmode\hskip 6mm plus 1mm minus 3mm\copy\sqwh\else\CASE
\hglue 6mm plus 1mm minus 3mm{\copy\sqwh\finalhyphendemerits=0\EMpar}\fi}
\def\cqfd{\ifmmode\hskip 6mm plus 1mm minus 3mm\copy\sqwh\else\CASE
\nobreak\hfil\penalty50\hskip1em\null\nobreak\hfil
{\copy\sqwh\parfillskip=0pt\finalhyphendemerits=0\let\par=\endgraf\EMpar}\fi}
%
%
\newcount\EMchap  \EMchap=0
\newcount\EMsec   \EMsec=0
\newcount\EMequa  \EMequa=0
\newcount\EMtheo  \EMtheo=0
\newread\infile   \newwrite\outfile
\def\latex{\openin\infile=\jobname.aux
\ifeof\infile\else\input\jobname.aux\fi
\immediate\openout\outfile=\jobname.aux}
\def\labchap#1{\global\EMsec=0  \global\EMequa=0 \global\EMtheo=0 
       \global\advance\EMchap by 1 {\the\EMchap}%
       \def\next{#1}{\ifx\next\empty\else\immediate\write\outfile%
       {\def\expandafter\noexpand\csname emEMeM#1\endcsname%
       {\the\EMchap}}\fi}}
\def\labsec#1{\global\advance\EMsec by 1 {\the\EMchap.\the\EMsec}%
       \def\next{#1}{\ifx\next\empty\else\immediate\write\outfile%
       {\def\expandafter\noexpand\csname emEMeM#1\endcsname%
       {\the\EMchap.\the\EMsec}}\fi}}
\def\labequa#1 {\global\advance\EMequa by 1 {\the\EMchap.\the\EMequa}%
       \def\next{#1}{\ifx\next\empty\else\immediate\write\outfile%
       {\def\expandafter\noexpand\csname emEMeM#1\endcsname%
       {\the\EMchap.\the\EMequa}}\fi}}
\def\labtheo#1{\global\advance\EMtheo by 1 {\the\EMchap.\the\EMtheo}%
       \def\next{#1}{\ifx\next\empty\else\immediate\write\outfile%
       {\def\expandafter\noexpand\csname emEMeM#1\endcsname%
       {\the\EMchap.\the\EMtheo}}\fi}}
\def\ref#1{{\csname emEMeM#1\endcsname}}
\let\?=\labequa
%
%
\def\EMdash{\hbox to 7.5mm{\vrule height .63ex depth -.59ex 
width 5.4mm\hfill}}
\newdimen\hangbiblio  \newcount\EMrefno  \EMrefno=1
\def\EMbegref#1#2{\removelastskip\goodbreak
\par\bgroup\EMrefno=1 \parindent=#2\petit\let\sl=\it
\vskip 10mm plus 6mm minus 3mm\goodbreak
\centerline{\maj{#1}}\nobreak\vskip\baselineskip\nobreak
\hangbiblio=\parindent \advance\hangbiblio by 7.5mm}
\def\EMrefa{\par\ifnum\EMrefno<3 \nobreak\else\goodbreak\fi
\vskip 1pt\noindent}
\def\EMrefb{\advance\EMrefno by 1 
\hangindent=\hangbiblio\hangafter=1 \ignorespaces }
%
%

%
%
\def\begrefnum#1{\EMbegref{#1}{10.5mm}\def\ref{\EMrefa
\hbox to\parindent{\hphantom{\,[00]}\llap{[\the\EMrefno]}\hfill}\EMrefb}%
\def\refsp##1.##2 {\EMrefa
\hbox to\parindent{\hphantom{\,[00}\llap{[##1}##2]\hfill}\EMrefb}}
%
%
\def\begreflab#1#2{\EMbegref{#1}{#2}\def\ref##1 {\EMrefa
\hbox to\parindent{[##1]\hfill}\EMrefb}}
\def\endref{\par\egroup}
%
%
\long\def\adresse#1{\removelastskip\nobreak\vskip\baselineskip\noindent
\vbox\bgroup\def\\{\vskip 1.0mm\parindent=.75cm\relax}\petit
\parindent=0pt\obeylines {#1}}
\def\endadresse{\egroup}
%
%
\def\begmat{\bgroup\parindent=6mm\rightskip=12.95mm\let\par=\endgraf
\def\s##1{\nobreak\item{##1}}%
\def\ss##1{\nobreak\itemitem{##1}}}
\def\pp#1 {\CASE\leaders\hbox to 2mm{\hfil.}\hfill
\rlap{\hbox to 12.95mm{\hfill #1}}\par}
%
%
\newcount\EMaaa  \newcount\EMbbb
\def\ssec#1#2{\removelastskip\goodbreak
\setbox0=\hbox{\smc\ignorespaces #1\CASE~}%
\dimen1=\wd0  \EMaaa=\dimen1  \dimen2=2mm  \EMbbb=\dimen2
\divide\EMaaa by \EMbbb  \advance\EMaaa by 1
\multiply\dimen2 by \EMaaa  \wd0=\dimen2
\medskip\noindent\hangindent=\dimen2\hangafter=1
\box0\ignorespaces #2}
\def\endmat{\par\egroup}
%
%
\def\begfig#1{\midinsert\kern #1}   
\def\endfig{\endinsert}
\newdimen\largeurfig  \largeurfig=\hsize
\newskip\deplacefig   \deplacefig=0mm      
\newskip\montefig     \montefig=0mm        
\newskip\EMpush
%
%
\def\figinsert#1{%
\EMpush=\hsize  \advance\EMpush by -\largeurfig
\divide\EMpush by 2  \advance\EMpush by \deplacefig
\par\nobreak\vskip-\montefig
\noindent\null\hskip\EMpush
\includegraphics{#1}
\vskip\montefig
}
%
%
\def\figure#1#2{\par\nobreak\vskip 5.0mm
\boitecentree{\pluspetit{\smc Figure}~\ignorespaces #1\CASE\par
\def\next{#2}\ifx\next\empty\else\nobreak\vskip 1.2mm 
\ignorespaces #2\CASE\fi}%
\vskip 1.2mm}
%
%
\def\doublefig#1#2#3{\bgroup\begfig{#1}\let\par=\endgraf
\line{\vtop{\hsize=.45\hsize #2}\hfill\vtop{\hsize=.45\hsize #3}}
\endfig\egroup}
\catcode`\@=\active
\newif\ifpremierepage  \premierepagetrue
\newtoks\headpremierepage  \headpremierepage={\hfil}
\newtoks\leftheadline  \newtoks\rightheadline  
\leftheadline={\tenpol\folio\hfil\bgroup\trespetit
\the\auteurcrt\egroup\hfil\hphantom{\folio}}
\rightheadline={\tenpol\hphantom{\folio}\hfil\bgroup\trespetit
\the\titrecrt\egroup\hfil\folio}
\headline={\ifpremierepage\rightheadline=\headpremierepage
\leftheadline=\headpremierepage\fi%
\ifodd\pageno\the\rightheadline%
\else\the\leftheadline\fi}
\footline={\ifpremierepage\global\premierepagefalse\fi%
\hfil\quad\hfil}
%
%
%
\font \tenrm                  = cmr10           
\font \ninerm                 = cmr9
\font \eightrm                = cmr8
\font \sevenrm                = cmr7
\font \sixrm                  = cmr6
\font \fiverm                 = cmr5
\font \teni                   = cmmi10          
\font \ninei                  = cmmi9
\font \eighti                 = cmmi8
\font \seveni                 = cmmi7
\font \sixi                   = cmmi6
\font \fivei                  = cmmi5
\font \tensy                  = cmsy10          
\font \ninesy                 = cmsy9
\font \eightsy                = cmsy8
\font \sevensy                = cmsy7
\font \sixsy                  = cmsy6
\font \fivesy                 = cmsy5
\font \tenex                  = cmex10          
\font \nineex                 = cmex9
\font \eightex                = cmex8
\font \sevenex                = cmex7
\font \sixex                  = cmex7 at 6pt
\font \fiveex                 = cmex7 at 5pt
\font \tenit                  = cmti10          
\font \nineit                 = cmti9
\font \eightit                = cmti8
\font \sevenit                = cmti7
\font \sixit                  = cmti7 at 6pt
\font \fiveit                 = cmti7 at 5pt
\font \tensl                  = cmsl10          
\font \ninesl                 = cmsl9
\font \eightsl                = cmsl8
\font \sevensl                = cmsl8 at 7pt
\font \sixsl                  = cmsl8 at 6pt
\font \fivesl                 = cmsl8 at 5pt
\font \tenbf                  = cmbx10          
\font \ninebf                 = cmbx9
\font \eightbf                = cmbx8
\font \sevenbf                = cmbx7
\font \sixbf                  = cmbx6
\font \fivebf                 = cmbx5
\font \tentt                  = cmtt10          
\font \ninett                 = cmtt9
\font \eighttt                = cmtt8
%
%
\font \tenmsa                 = msam10          
\font \ninemsa                = msam9
\font \eightmsa               = msam8
\font \sevenmsa               = msam7
\font \sixmsa                 = msam6
\font \fivemsa                = msam5
\font \tenmsb                 = msbm10          
\font \ninemsb                = msbm9
\font \eightmsb               = msbm8
\font \sevenmsb               = msbm7
\font \sixmsb                 = msbm6
\font \fivemsb                = msbm5
%
%
\newfam\eufmfam                                 
\newfam\sansfam                                 
\font \teneufm                = eufm10
\font \nineeufm               = eufm9
\font \eighteufm              = eufm8
\font \seveneufm              = eufm7
\font \sixeufm                = eufm6
\font \fiveeufm               = eufm5
\font \tensans                = cmss10            
\font \ninesans               = cmss9
\font \eightsans              = cmss8
\font \sevensans              = cmss8 at 7pt
\font \sixsans                = cmss8 at 6pt
\font \fivesans               = cmss8 at 5pt
\font \tensmc                 = cmcsc10
\font \ninesmc                = cmcsc9
\font \eightsmc               = cmcsc8

\font \tenbfsl                = cmbxsl10           
\font \ninebfsl               = cmbxsl10 at 9pt
\font \eightbfsl              = cmbxsl10 at 8pt

\skewchar\ninei='177 \skewchar\eighti='177 \skewchar\sixi='177
\skewchar\ninesy='60 \skewchar\eightsy='60 \skewchar\sixsy='60
\hyphenchar\ninett=-1
\def\loadeuex{\Loadmathfont{euex}\let\EUEX=Y}
\def\loadeurm{\Loadmathfont{eurm}\let\EURM=Y}
\def\loadeusm{\Loadmathfont{eusm}\def\scr##1{{\eusm ##1}}\let\EUSM=Y}
\def\loadeufb{\Loadmathfont{eufb}\let\EUFB=Y}
\def\loadeurb{\Loadmathfont{eurb}\let\EURB=Y}
\def\loadeusb{\Loadmathfont{eusb}\let\EUSB=Y}
\catcode`\@=11
\def\Loadmathfont#1{%
   \expandafter\font\csname ten#1\endcsname=#110
   \expandafter\font\csname nine#1\endcsname=#19
   \expandafter\font\csname eight#1\endcsname=#18
   \expandafter\font\csname seven#1\endcsname=#17
   \expandafter\font\csname six#1\endcsname=#16
   \expandafter\font\csname five#1\endcsname=#15
   \edef\next{\noexpand\alloc@@8\fam\chardef\sixt@@n
     \expandafter\noexpand\csname#1fam\endcsname}%
   \next
   \textfont\csname#1fam\endcsname \csname ten#1\endcsname
   \scriptfont\csname#1fam\endcsname \csname seven#1\endcsname
   \scriptscriptfont\csname#1fam\endcsname \csname five#1\endcsname
   \expandafter\def\csname #1\expandafter\endcsname\expandafter{%
      \expandafter\fam\csname#1fam\endcsname
      \csname ten#1\endcsname}%
   \expandafter\gdef\csname load#1\endcsname{}%
}%
\catcode`\@=\active
\def\famsupp{%
\ifx\EUEX\undefined\else\textfont\euexfam=\teneuex 
\scriptfont\euexfam=\seveneuex \scriptscriptfont\euexfam=\fiveeuex 
\def\euex{\fam\euexfam\teneuex}\fi
\ifx\EURM\undefined\else\textfont\eurmfam=\teneurm 
\scriptfont\eurmfam=\seveneurm \scriptscriptfont\eurmfam=\fiveeurm 
\def\eurm{\fam\eurmfam\teneurm}\fi
\ifx\EUSM\undefined\else\textfont\eusmfam=\teneusm 
\scriptfont\eusmfam=\seveneusm \scriptscriptfont\eusmfam=\fiveeusm 
\def\eusm{\fam\eusmfam\teneusm}\fi
\ifx\EUFB\undefined\else\textfont\eufbfam=\teneufb 
\scriptfont\eufbfam=\seveneufb \scriptscriptfont\eufbfam=\fiveeufb 
\def\eufb{\fam\eufbfam\teneufb}\fi
\ifx\EURB\undefined\else\textfont\eurbfam=\teneurb 
\scriptfont\eurbfam=\seveneurb \scriptscriptfont\eurbfam=\fiveeurb 
\def\eurb{\fam\eurbfam\teneurb}\fi
\ifx\EUSB\undefined\else\textfont\eusbfam=\teneusb 
\scriptfont\eusbfam=\seveneusb \scriptscriptfont\eusbfam=\fiveeusb 
\def\eusb{\fam\eusbfam\teneusb}\fi}
\def\famsupppetit{%
\ifx\EUEX\undefined\else\textfont\euexfam=\nineeuex 
\scriptfont\euexfam=\sixeuex \scriptscriptfont\euexfam=\fiveeuex 
\def\euex{\fam\euexfam\nineeuex}\fi
\ifx\EURM\undefined\else\textfont\eurmfam=\nineeurm 
\scriptfont\eurmfam=\sixeurm \scriptscriptfont\eurmfam=\fiveeurm 
\def\eurm{\fam\eurmfam\nineeurm}\fi
\ifx\EUSM\undefined\else\textfont\eusmfam=\nineeusm 
\scriptfont\eusmfam=\sixeusm \scriptscriptfont\eusmfam=\fiveeusm 
\def\eusm{\fam\eusmfam\nineeusm}\fi
\ifx\EUFB\undefined\else\textfont\eufbfam=\nineeufb 
\scriptfont\eufbfam=\sixeufb \scriptscriptfont\eufbfam=\fiveeufb 
\def\eufb{\fam\eufbfam\nineeufb}\fi
\ifx\EURB\undefined\else\textfont\eurbfam=\nineeurb 
\scriptfont\eurbfam=\sixeurb \scriptscriptfont\eurbfam=\fiveeurb 
\def\eurb{\fam\eurbfam\nineeurb}\fi
\ifx\EUSB\undefined\else\textfont\eusbfam=\nineeusb 
\scriptfont\eusbfam=\sixeusb \scriptscriptfont\eusbfam=\fiveeusb 
\def\eusb{\fam\eusbfam\nineeusb}\fi}
\def\famsupppluspetit{%
\ifx\EUEX\undefined\else\textfont\euexfam=\eighteuex 
\scriptfont\euexfam=\sixeuex \scriptscriptfont\euexfam=\fiveeuex 
\def\euex{\fam\euexfam\eighteuex}\fi
\ifx\EURM\undefined\else\textfont\eurmfam=\eighteurm 
\scriptfont\eurmfam=\sixeurm \scriptscriptfont\eurmfam=\fiveeurm 
\def\eurm{\fam\eurmfam\eighteurm}\fi
\ifx\EUSM\undefined\else\textfont\eusmfam=\eighteusm 
\scriptfont\eusmfam=\sixeusm \scriptscriptfont\eusmfam=\fiveeusm 
\def\eusm{\fam\eusmfam\eighteusm}\fi
\ifx\EUFB\undefined\else\textfont\eufbfam=\eighteufb 
\scriptfont\eufbfam=\sixeufb \scriptscriptfont\eufbfam=\fiveeufb 
\def\eufb{\fam\eufbfam\eighteufb}\fi
\ifx\EURB\undefined\else\textfont\eurbfam=\eighteurb 
\scriptfont\eurbfam=\sixeurb \scriptscriptfont\eurbfam=\fiveeurb 
\def\eurb{\fam\eurbfam\eighteurb}\fi
\ifx\EUSB\undefined\else\textfont\eusbfam=\eighteusb 
\scriptfont\eusbfam=\sixeusb \scriptscriptfont\eusbfam=\fiveeusb 
\def\eusb{\fam\eusbfam\eighteusb}\fi}
%
%

%
%

\def\desserre{\mathsurround=.1em \relax}

\def\desserreplus{\mathsurround=.115em \relax}
\font \tenpol                 = cmr10
\font \ninepol                = cmr9

\font \sixpol                 = cmr6
%

%

%
%
\catcode`\@=11
\def\normal{%
\textfont0=\tenrm \scriptfont0=\sevenrm \scriptscriptfont0=\fiverm
\def\rm{\fam0\tenrm}%
\textfont1=\teni \scriptfont1=\seveni \scriptscriptfont1=\fivei
\def\oldstyle{\fam1\teni}%
\textfont2=\tensy \scriptfont2=\sevensy \scriptscriptfont2=\fivesy
\textfont3=\tenex \scriptfont3=\sevenex \scriptscriptfont3=\fiveex
\textfont\itfam=\tenit \scriptfont\itfam=\sevenit
\scriptscriptfont\itfam=\fiveit \def\it{\fam\itfam\tenit}%
\textfont\slfam=\tensl \scriptfont\slfam=\sevensl
\scriptscriptfont\slfam=\fivesl \def\sl{\fam\slfam\tensl}%
\textfont\bffam=\tenbf \scriptfont\bffam=\sevenbf
\scriptscriptfont\bffam=\fivebf \def\bf{\fam\bffam\tenbf}%
\textfont\ttfam=\tentt \def\tt{\fam\ttfam\tentt}%
\textfont\msafam=\tenmsa \scriptfont\msafam=\sevenmsa
\scriptscriptfont\msafam=\fivemsa
\textfont\msbfam=\tenmsb \scriptfont\msbfam=\sevenmsb
\scriptscriptfont\msbfam=\fivemsb
\textfont\eufmfam=\teneufm \scriptfont\eufmfam=\seveneufm
\scriptscriptfont\eufmfam=\fiveeufm \def\goth{\fam\eufmfam\teneufm}%
\textfont\sansfam=\tensans \scriptfont\sansfam=\sevensans
\scriptscriptfont\sansfam=\fivesans \def\sans{\fam\sansfam\tensans}%
\famsupp
\def\smc{\tensmc}%
\def\bfsl{\tenbfsl}%
\normalbaselineskip=4.5mm plus .02mm minus .01mm
\EMskip=\normalbaselineskip  \advance\EMskip by 0mm plus .5mm minus .1mm
\setbox\strutbox=\hbox{\vrule height3.2mm depth1.3mm width0pt}%
\def\big##1{{\hbox{$\left##1\vbox to8.5\p@{}\right.\n@space$}}}%
\def\Big##1{{\hbox{$\left##1\vbox to11.5\p@{}\right.\n@space$}}}%
\def\bigg##1{{\hbox{$\left##1\vbox to14.5\p@{}\right.\n@space$}}}%
\def\Bigg##1{{\hbox{$\left##1\vbox to17.5\p@{}\right.\n@space$}}}%
\normalbaselines\rm\makesqwh
}%
\def\petit{%
\textfont0=\ninerm \scriptfont0=\sixrm \scriptscriptfont0=\fiverm
\def\rm{\fam0\ninerm}%
\textfont1=\ninei \scriptfont1=\sixi \scriptscriptfont1=\fivei
\def\oldstyle{\fam1\ninei}%
\textfont2=\ninesy \scriptfont2=\sixsy \scriptscriptfont2=\fivesy
\textfont3=\nineex \scriptfont3=\sixex \scriptscriptfont3=\fiveex
\textfont\itfam=\nineit \scriptfont\itfam=\sixit
\scriptscriptfont\itfam=\fiveit \def\it{\fam\itfam\nineit}%
\textfont\slfam=\ninesl \scriptfont\slfam=\sixsl
\scriptscriptfont\slfam=\fivesl \def\sl{\fam\slfam\ninesl}%
\textfont\bffam=\ninebf \scriptfont\bffam=\sixbf
\scriptscriptfont\bffam=\fivebf \def\bf{\fam\bffam\ninebf}%
\textfont\ttfam=\ninett \def\tt{\fam\ttfam\ninett}%
\textfont\msafam=\ninemsa \scriptfont\msafam=\sixmsa
\scriptscriptfont\msafam=\fivemsa
\textfont\msbfam=\ninemsb \scriptfont\msbfam=\sixmsb
\scriptscriptfont\msbfam=\fivemsb
\textfont\eufmfam=\nineeufm \scriptfont\eufmfam=\sixeufm
\scriptscriptfont\eufmfam=\fiveeufm \def\goth{\fam\eufmfam\nineeufm}%
\textfont\sansfam=\ninesans\scriptfont\sansfam=\sixsans
\scriptscriptfont\sansfam=\fivesans \def\sans{\fam\sansfam\ninesans}%
\famsupppetit
\def\smc{\ninesmc}%
\def\bfsl{\ninebfsl}%
\normalbaselineskip=3.5mm plus .015mm
\EMskip=\normalbaselineskip  \advance\EMskip by 0mm plus .4mm
\setbox\strutbox=\hbox{\vrule height2.5mm depth1.0mm width0pt}%
\def\big##1{{\hbox{$\left##1\vbox to7.65\p@{}\right.\n@space$}}}%
\def\Big##1{{\hbox{$\left##1\vbox to10.35\p@{}\right.\n@space$}}}%
\def\bigg##1{{\hbox{$\left##1\vbox to13.05\p@{}\right.\n@space$}}}%
\def\Bigg##1{{\hbox{$\left##1\vbox to15.75\p@{}\right.\n@space$}}}%
\normalbaselines\rm\makesqwh
}%
\def\pluspetit{%
\textfont0=\eightrm \scriptfont0=\sixrm \scriptscriptfont0=\fiverm
\def\rm{\fam0\eightrm}%
\textfont1=\eighti \scriptfont1=\sixi \scriptscriptfont1=\fivei
\def\oldstyle{\fam1\eighti}%
\textfont2=\eightsy \scriptfont2=\sixsy \scriptscriptfont2=\fivesy
\textfont3=\eightex \scriptfont3=\sixex \scriptscriptfont3=\fiveex
\textfont\itfam=\eightit \scriptfont\itfam=\sixit
\scriptscriptfont\itfam=\fiveit \def\it{\fam\itfam\eightit}%
\textfont\slfam=\eightsl \scriptfont\slfam=\sixsl
\scriptscriptfont\slfam=\fivesl \def\sl{\fam\slfam\eightsl}%
\textfont\bffam=\eightbf \scriptfont\bffam=\sixbf
\scriptscriptfont\bffam=\fivebf \def\bf{\fam\bffam\eightbf}%
\textfont\ttfam=\eighttt \def\tt{\fam\ttfam\eighttt}%
\textfont\msafam=\eightmsa \scriptfont\msafam=\sixmsa
\scriptscriptfont\msafam=\fivemsa
\textfont\msbfam=\eightmsb \scriptfont\msbfam=\sixmsb
\scriptscriptfont\msbfam=\fivemsb
\textfont\eufmfam=\eighteufm \scriptfont\eufmfam=\sixeufm
\scriptscriptfont\eufmfam=\fiveeufm \def\goth{\fam\eufmfam\eighteufm}%
\textfont\sansfam=\eightsans\scriptfont\sansfam=\sixsans
\scriptscriptfont\sansfam=\fivesans \def\sans{\fam\sansfam\eightsans}%
\famsupppluspetit
\def\smc{\eightsmc}%
\def\bfsl{\eightbfsl}%
\normalbaselineskip=3mm plus .01mm
\EMskip=\normalbaselineskip  \advance\EMskip by 0mm plus .3mm
\setbox\strutbox=\hbox{\vrule height2.1mm depth.9mm width0pt}%
\def\big##1{{\hbox{$\left##1\vbox to6.8\p@{}\right.\n@space$}}}%
\def\Big##1{{\hbox{$\left##1\vbox to9.2\p@{}\right.\n@space$}}}%
\def\bigg##1{{\hbox{$\left##1\vbox to11.6\p@{}\right.\n@space$}}}%
\def\Bigg##1{{\hbox{$\left##1\vbox to14.0\p@{}\right.\n@space$}}}%
\normalbaselines\rm\makesqwh
}%
\def\trespetit{\pluspetit}%
%
%
\def\frak{\relaxnext@\ifmmode\let\next\frak@\else
 \def\next{\Err@{Use \string\frak\space only in math mode}}\fi\next}
\def\frak@#1{{\frak@@{#1}}}
\def\frak@@#1{\noaccents@\fam\eufmfam#1}
\catcode`\@=\active
%
%
\let\EMsavS=\S
\def\S{\leavevmode\unkern\EMsavS\kern.1em\ignorespaces}
\def\og{\leavevmode\raise.3ex\hbox{$\scriptscriptstyle\langle\!\langle$}}
\def\fg{\leavevmode\raise.3ex\hbox{$\scriptscriptstyle\rangle\!\rangle$}}

\normal
%

%
%
%
%
\ifx\MACUTIL\undefined\let\MACUTIL=Y\else  \fi
%
%

%
%

%
%

%
%

%
%

%
%

%
%

%
%
%
%

%
%
\def\norm#1{{\left\Vert{@,@, #1 @,@,}\right\Vert}}      
%
%

%
%
\def\Aut{\operatorname{Aut}}



\def\SU{\operatorname{SU}}

%
%


\def\up#1{\leavevmode\raise.18ex\hbox{#1}}

\def\longrarr{\hbox to 15mm{\enspace\rightarrowfill}}
\def\longdarr{\hbox to 15mm{\hfill$\relbar\relbar\to$\hfill}}
\def\darr{\mathrel{\downarrow}}
\setbox111=\hbox{\raise.65ex\hbox to 5mm{\hfill$\shortmid$\hfill}}
\setbox112=\hbox{\raise1.2ex\hbox to 5mm{\hfill$\shortmid$\hfill}}
\setbox113=\hbox to 5mm{\hfill$\darr$\hfill}

%
%
\hyphenation{Sprin-ger
ab-solu ab-so-lue ab-so-lus ab-so-lues ab-so-lu-ment
as-so-cia-tion as-so-cia-tions
aupa-ra-vant auto-bio-gra-phi-que auto-bio-gra-phi-ques
bio-logie bio-lo-gique bio-lo-gi-ques
civi-li-sa-tion civi-li-sa-tions civi-li-ser
com-mer-cial com-mer-ciale com-mer-ciaux com-mer-ciales
con-ve-nu con-ve-nue con-ve-nus con-ve-nues
cor-res-pon-dance
di-mi-nuer en-cou-rait en-cou-raient
exa-mi-nera exem-ple exem-ples exis-te exis-ter
ima-gi-ne ima-gi-na-tion ima-gi-naire im-por-tance
in-dien in-diens in-dien-ne in-dien-nes
in-tui-tif in-tui-tifs in-tui-tive in-tui-tives in-tui-ti-ve-ment
jeune- jeu-nes
la-quel-le li-ber-taire li-ber-tai-res long-temps
mani-pu-la-tion mani-pu-la-tions mo-derne mo-der-nes mo-der-ni-sa-tion
mo-di-fier mo-di-fie mo-di-fient
mon-ta-gne mon-ta-gnes mon-ta-gneux mon-ta-gneu-se mon-ta-gneu-ses
mytho-lo-gique mytho-lo-giques
orien-ta-lisme orien-ta-lismes
pater-na-lisme per-met-tre per-met-tent per-met-tant
po-pu-laire po-pu-lai-res po-pu-la-tion po-pu-la-tions
pri-son-nier pri-son-niers puis-que quel-que quel-ques
recher-che recher-ches recon-nais-sable recon-nais-sa-bles
re-de-va-ble re-de-va-bles
re-la-tif re-la-tifs re-la-tive re-la-tives re-la-ti-ve-ment
res-pon-sable res-pon-sa-bles
sa-tis-fait sa-tis-faite sa-tis-faits sa-tis-fai-tes sa-tis-fai-re
se-con-daire se-con-dai-res
si-tua-tion si-tua-tions suf-fi-sam-ment sui-vant sui-vent sui-vre
tota-li-taire tota-li-taires ty-pi-que ty-pi-ques ty-pi-que-ment
vien-nent vio-len-ce voca-bu-laire voca-bu-lai-res
}
%
%
\def\car.{carac\-t\'e\-ris\-tique}
\def\dev.{d\'e\-ve\-lop\-p}
\def\extr.{ex\-t\'e\-rieur}
\def\intr.{in\-t\'e\-rieur}
\let\phi=\varphi
 
\def\abs#1{\vert#1\vert}
\def\norm#1{\left\Vert#1\right\Vert}
\def\R {{\Bbb R}}
\def\Q {{\Bbb Q}}
\def\I {{\Bbb I}}
\def\C {{\Bbb C}}
\def\N{{\Bbb N}}
\def\Ur{{\Bbb U}}
\def\e{{\varepsilon}}

\def\Z {{\Bbb Z}}

\def\SU{{\operatorname{SU}}}
\def\SO{{\operatorname{SO}}}
\def\HL{{\underline{H}_1{\Cal L}\iota}}
\def\U{{\operatorname{U}}}

\def\Iso{{\operatorname{Iso}}\,}
\def\Aut{{\operatorname{Aut}}}

\def\LO{{\operatorname{LO}}}

\def\UCB{{\operatorname{UCB}}\,}

\def\UCB{{\operatorname{UCB}}}

\def\H{{\Cal H}}

\def\s{{\Bbb S}}


\latex

\titre{\maj{$mm$-Spaces and group actions}}
 
\auteurbref{V. Pestov}
 
\auteur{Vladimir \sn{Pestov}}
 
\abstract{These are introductory notes on some aspects of concentration of
measure in the presence of an acting group and its links
to Ramsey 
theory\note{
Based on a
lecture given in the framework of
{\it
S\'eminaire Borel de IIIe Cycle romand de Math\'ematiques: 
``2001: an mm-space odyssey''
{\rm (}Espaces avec une m\'etrique et une mesure, d'apr\`es 
M. Gromov{\rm )}} at the Institute of Mathematics,
University of Bern and
a {\it S\'eminaire du Li\`evre} talk at the
Department of Mathematics, University of Geneva.
The author gratefully acknowledges generous
support from the Swiss National Science Foundation during his visit
in April--May 2001 and thanks
Pierre de la Harpe for his hospitality and many stimulating
conversations. While in Switzerland, the 
author has also greatly benefitted from discussions
with Gulnara Arzhantseva, Anna Erschler,
Thierry Giordano, Eli Glasner, 
Rostislav Grigorchuk, Volodymyr Nekrashevych, Vitali Milman,
and Tatiana Nagnibeda.
Partial support also came from the Marsden Fund of
the Royal Society of New Zealand. Numerous remarks by
the anonymous referee have been most helpful.}.}

\sec{\\/}{Introduction}

It can be argued that the theory we are
interested in (call it theory of $mm$-spaces, the phenomenon of 
concentration of measure
on high-dimensional structures, asymptotic
geometric analysis, geometry of large dimensions$\ldots$) has been 
largely shaped up by three publications.
These are: the book by Paul L\'evy
\cite{L\'ev}, Vitali Milman's new proof of the
Dvoretzky theorem \cite{M1}, and the paper by Gromov and Milman
\cite{Gr-M1} which had set up a framework for systematically dealing
with concentration of measure. 
Significantly, in the two latter papers concentration
goes hand in hand with group actions on suitable spaces with
metric and measure.

It is also known that concentration of measure
and combinatorial, Ramsey-type results have a similar nature
and are often found together \cite{M3}.

A number of attempts have been made to understand the nature of
the interplay between concentration, transformation groups, 
and Ramsey theory,
cf. papers by Gromov \cite{Gr1}, Milman
\cite{M2,M3}, and some others \cite{A-M,Gl,P2,P3,G-P,Gl-W}.
However, it is safe to say that there is still a long way to go
towards the full understanding of the picture.

Here we aim at providing a readable introduction into this circle
of ideas.

\sec{\\/}{Some concepts of asymptotic geometric analysis}

\parag{Definition~1}
A {\it space with metric and measure,} or an $mm${\it -space}, is a
triple $(X,d,\mu)$, where $d$ is a metric on a set $X$ and $\mu$ is
a finite Borel measure on the metric space $(X,d)$. It will be
convenient to assume throughout that $\mu$ is a
probability measure, that is, normalized to one. 

\parag{Definition~2}
The {\it concentration function} $\alpha_X$ 
of an $mm$-space $X=(X,d,\mu)$ is
defined for non-negative real $\e$ as follows:
$$
\alpha_X(\e)=\cases\frac 12, & \text{if $\e=0$,} \\
1-\inf\{\mu(A_\e)\colon \text{$A\subseteq X$ is Borel, $\mu(A)\geq\frac 12$}\},
&\text{if $\e>0$.}
\endcases
$$
Here by $A_\e$ we denote the $\e$-neighbourhood ($\e$-fattening,
$\e$-thickening) of $A$.

\parag{Exercise~1}
Prove that $\alpha(\e)\to 0$ as $\e\to\infty$. (For spaces of finite
diameter this is of course obvious.)

\parag{Definition~3}
An infinite family of $mm$-spaces, $(X_n,d_n,\mu_n)_{n=1}^\infty$,
is called a {\it L\'evy family} if the concentration functions $\alpha_n$ of
$X_n$ converge to zero pointwise on $(0,\infty)$:
$$
\forall \e>0,~~ \alpha_n(\e)\to 0 \text{ as } n\to\infty.
$$

\parag{Exercise~2}
Prove that the above condition is equivalent to the following.
Let $A_n\subseteq X_n$ be Borel subsets with the property that
$$
\liminf_{n\to\infty}\mu_n(A_n)>0.
$$
Then
$$
\forall \e>0,~~\lim_{n\to\infty}\mu_n((A_n)_\e)=1.
$$
\vskip .3cm

The following are some of
the most common examples of L\'evy families.

\parag{Example 1}
Unit spheres $\s^n$ in the Euclidean spaces $\R^{n+1}$, equipped with the
Euclidean (or geodesic) distances and the
normalized
Haar measures (that is, the unique rotation-invariant probability measures).
This result is due to Paul L\'evy \cite{L\'ev}, though his proof, based on
the isoperimetric inequality,
was only made rigorous much later by Gromov \cite{Gr2}. 
(Nowadays simpler proofs, using the Brunn--Minkowski 
inequality, are known, cf. \cite{Gr-M2, Sch}.)

\parag{Example 2}
The special orthogonal groups $\SO(n)$, equipped with the normalized
Haar measure and the uniform operator metric,
$$
d(T,S):=\norm{T-S},
$$
induced from ${\Cal B}(\R^n)\cong M_n$. 
This was established by Gromov and Milman \cite{Gr-M1}.
The same argument holds for the special unitary groups.


\parag{Example~3}
The family of finite permutation groups $(S_n)$, 
equipped with the uniform (normalized counting)
measure and the Hamming distance:
$$
d(\sigma,\tau)=\frac 1n \vert\{i\colon \sigma(i)\neq \tau(i)\}\vert.
$$
The result is due to Maurey \cite{Ma}, see also \cite{Ta1}.

\parag{Example~4} The Hamming cubes $\{0,1\}^n$ equipped with
the normalized counting measure and the Hamming distance 
$d(x,y)=\frac 1n\vert\{i\colon x_i\neq y_i\}\vert$ form
a L\'evy family \cite{Sch,M-S}.

\parag{Remark~1} All of the above are
{\it normal} L\'evy families, meaning that the concentration
functions $\alpha_n$ admit Gaussian upper bounds:
$$
\alpha_n(\e)\leq C_1\exp(-C_2n\e^2)
$$
for some $C_1,C_2>0$.

It should be noted that this is not always the case for `naturally occuring'
L\'evy families. For instance, the groups 
$\operatorname{SL}(2,{\Bbb F}_p)$, 
where $p$ are prime numbers, 
equipped with the normalized counting
measure and the word metric given by a fixed system of generators
in $\operatorname{SL}(2,{\Bbb Z})$, form a L\'evy family
with $\alpha_p(\e)\leq C_1 \exp(-C_2\sqrt p\,\e)$,
\cite{A-M, M4}. (Recall in this connection that 
the $n$-th prime number $p_n\sim n\log n$.)

\parag{Remark~2} In Example 4, replace $\{0,1\}$ with
any probability measure space, $X=(X,\mu)$.
Equip every finite
power $X^n$ with the product measure $\mu^{\otimes\,n}$
and the normalized Hamming distance 
$d(x,y)=\frac 1n\vert\{i\colon x_i\neq y_i\}\vert$.
Unless $X$ is purely atomic,
the measures $\mu^{\otimes\,n}$ are not Borel,
and thus $X^n$ aren't even $mm$-spaces in the
sense of our definition. At the same time, if in the definition of
the concentration function 
we only restrict ourselves
to measurable subsets $A$ such that $A_\e$ are also measurable,
it can be shown
that $X^n,n\in\N$ form a L\'evy family in a very reasonable sense. 
(See \cite{Ta1,Ta3} for
far-reaching variations.) If anything, this shows that the
full formalization of the subject has not yet been achieved and nothing
is cast in stone.
\vskip .3cm

Notice that the $mm$-spaces from the above examples 1--4
are at the same time (phase spaces of) 
topological transformation groups,
with both metrics and measures being invariant
under group actions. In example 1 it is the action of the
orthogonal --- or the unitary --- group on the sphere, while in examples
2--4 the groups act upon themselves on the left.

\sec{\\/}{A transformation group framework}
Here is the idea of what kind of interaction between 
concentration phenomenon and
group actions one should expect. The following example is borrowed from
a paper by Vitali Milman \cite{M4}.

Suppose a group $G$ acts on an $mm$-space $(X,d,\mu)$ by
measure-preserving isometries.
Assume that the $mm$-space $X$ strongly concentrates, that is, 
the function
$\alpha_X(\e)$ drops off sharply already for small values of $\e$.
Let us assume, for instance, that
the concentration is so strong that, whenever $\mu(A)\geq {\frac {1}{7}}$, 
the measure of
the ${\frac {1}{10}}$-neighbourhood of $A$ is strictly
greater than $0.\,99$. 
(Cf. Exercise 2.) 

If now we partition $X$ into seven pieces, and pick at random
one hundred elements $g_1,g_2,\dots,g_{100}\in G$,
then at least 
one of the pieces, say $A$, has the property that all one hundred translates,
of ${\frac {1}{10}}$-neighbourhoods of $A$ by our elements $g_i$
have a point, $x^\ast$, in common. 
Equivalently, $x^\ast$ is `close' (closer than
${\frac {1}{10}}$) to each of the one hundred translates of $A$. 

The above effect becomes more
pronounced the higher the level of concentration is. Partition a
concentrated (`high-dimensional')
$mm$-space into a small number of subsets, and
at least one of them is hard to move.

\begfig{3cm}
\largeurfig=5cm
\deplacefig=1mm
\montefig=-1mm
\figinsert{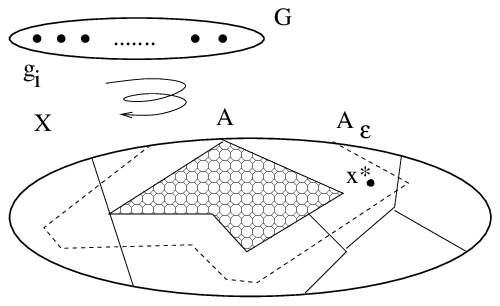}
\figure1{Dynamics in the presence of concentration}
\endfig

In order to set up a formal framework,
we assume all topological spaces and topological groups appearing in
this article to be metrizable, for the reasons of mere
technical simplicity.\note{More generally, metrics can be replaced
with uniform structures.}
We need $G$-spaces of a particular kind.
Let $X=(X,d)$ be a metric space,
not necessarily compact, and let
a group $G$ act on $X$ (on the left) by {\it uniformly} continuous maps.
In other words, there is
a map $G\times X\to X$, $(g,x)\mapsto g\cdot x$, such that
$g\cdot (h\cdot x) = (gh)\cdot x$, $e\cdot x =x $, and every map of the form
$$
X\ni x\mapsto g\cdot x \in X
$$
(a translation by $g$) is uniformly continuous. 
(Then it is automatically a uniform isomorphism.)  
If, moreover, $G$ is a topological group, then we require the action
$G\times X\to X$ to be continuous. 

\parag{Example~5} 
The motivation for our choice of
the class of $G$-spaces is provided by
the fact that every (metrizable)
compact $G$-space, $K$, is such: a translation
of $K$ by an element $g\in G$, being
a continuous map on a compact space, is uniformly continuous.
\vskip .3cm

Here is another property that compact $G$-spaces
possess automatically, while $G$-spaces of a more general nature do not.

\parag{Exercise~3} Let a topological group $G$ act continuously on a 
(metrizable) compact space $K=(K,d)$. Prove that for every $\e>0$
there is a neighbourhood of identity $V\ni e_G$ with the property that
whenever $g\in V$ and $x\in K$, one has $d(x,g\cdot x)<\e$.
[In abstract topological dynamics such actions are termed {\it bounded},
or else {\it motion equicontinuous}.]
\par
[{\it Hint:} using the continuity of the action $G\times K\to K$, 
choose for each
$x\in K$ a neighbourhood $U_x$ of $x$ in $K$ and a neighbourhood, $V_x$,
of $e_G$ in $G$, such that $V_x\cdot U_x\subseteq B_{\e}(x)$
(the open $d$-ball around $x$); now select a finite subcover of
$\{U_x\}\ldots$]

\parag{Example~6} 
Every metrizable group admits a right-invariant compatible metric
($d(x,y)=d(xa,ya)$), as well as a left-invariant one
($d(x,y)=d(ax,ay)$).
The action of $G$ on itself by 
left translations is an action by isometries with respect to a
left-invariant metric, and (exercise) 
an action by uniform isomorphisms with respect
to a right-invariant metric. 

\parag{Exercise~4}
Show that the action of a topological group $G$ upon itself, equipped
with a right invariant metric,  by left
translations, is bounded.

\parag{Example~7}
One topological group of interest to us is $\U(\H)_s$,
the full unitary group of a separable Hilbert space with the
{\it strong} operator topology. (That is, the topology induced from the 
Tychonoff product $\H^{\H}$.) A standard neighbourhood of identity in this
topology consists of all $T\in \U(\H)$ such that $\norm{T(x_i)-x_i}<\e$
for $i=1,2,\dots,n$, where $x_1,\dots,x_n$ is a finite 
collection of unit vectors in $\H$. This topology on $\U(\H)$
 coincides with
the {\it weak operator topology,} that is, the weakest topology
making continuous every map of the form
$$
\U(\H)\ni T\mapsto \langle x,Tx\rangle\in\C,~~ x\in\H.
$$

\parag{Example~8}
Let $\pi$ be a unitary representation of a group $G$
(viewed as discrete) in a Hilbert space
$\H$. Denote by $\s^{\infty}$ the unit sphere in $\H$, equipped with the norm
distance. Then $G$ acts on $\s^{\infty}$ by
isometries: $(g,x)\mapsto \pi_gx$. 

\parag{Remark~3} The above $G$-space is bounded for trivial reasons. It should
be noted, however, that in general one does not expect a `typical' $G$-space to
be bounded at all.

\parag{Definition~4}
Let a topological group $G$ act continuously, by uniform isomorphisms,
on two metric spaces, $X$ and $Y$. A {\it morphism}, or an 
{\it equivariant map,} from $X$ to $Y$ is a uniformly continuous
map $i\colon X\to Y$ which commutes with the action:
$$
i(g\cdot x) = g\cdot i(x).
$$

\parag{Definition~5} Let a topological group $G$ act continuously on 
a metric space $(X,d)$ by uniformly continuous maps, and let also $G$ act
continuously on a compact space $K$. 
Let  $i\colon X\to K$ be a morphism of $G$-spaces with an
everywhere dense image in $K$. The pair $(K,i)$ is called an 
{\it equivariant compactification} of $X$.

\parag{Example~9} Let $G$ and $\H$ be as in Example 8.
The unit ball $\Bbb B$ in $\H$ equipped with the weak
topology is compact, and $G$ acts on $\Bbb B$ in the same way as on the
sphere. The embedding $\s^{\infty}\hookrightarrow {\Bbb B}$ is
an equivariant compactification. 
\vskip .2cm
The following is at the heart of abstract
topological dynamics.

\theorem{1}{Let $G$ be a topological group, and let $d$ be a
right-invariant metric generating the topology of $G$. 
Let $K$ be a {\rm(}metric{\rm)} compact $G$-space, and let $\kappa\in K$ be
arbitrary. There is a morphism of $G$-spaces $i\colon (G,d)\to K$
such that $i(e)=\kappa$.}

\proof{
Define the map $i\colon G\to K$ (an {\it orbit map}) by
$$i\colon G\ni y\mapsto y\cdot \kappa\in K.$$
This map is equivariant. [$i(g\cdot y) = (gy)\cdot \kappa
=g\cdot (y\cdot \kappa)= g\cdot i(y)$.]
It only remains to check the uniform continuity of $i$. Choose any
continuous metric on $K$, say $\rho$. Using Exercise 3,
find a $\delta>0$ with the property that
$\rho(x,g\cdot x)<\e$ whenever $x\in K$ and $d(g,e_G)<\delta$. 
If now $g,h\in G$ are such that $d(g,h)<\delta$, then $d(gh^{-1},e_G)<\delta$
and consequently
$$\rho(h\kappa,g\kappa) =\rho(h\kappa,gh^{-1}(h\kappa))
<\e.$$
\qed}

\parag{Remark~4} 
The difference between the right and left invariant metrics (or, more
generally, uniform structures) on a topological group cannot be 
overemphasized. Even if they are totally symmetric, they cease to be such
as soon as we choose the action (in our case, by left
translations).
\parag{}
Here is a key notion putting the concentration of measure in a
dynamical context.

\parag{Definition~6}
Let a metrizable
topological group $G$ act continuously
by uniform isomorpisms on a metric space $X=(X,d)$. 
Say that the $G$-space (transformation group)
$(G,X)$ is {\it L\'evy} (Gromov and Milman
\cite{Gr-M1}) if
there are a sequence of subgroups of $G$
$$
G_1\subseteq G_2\subseteq \dots \subseteq G_n \subseteq \dots \subseteq G,
$$
and a sequence of probability measures 
$$
\mu_1,\mu_2,\dots,\mu_n,\dots
$$
on $(X,d)$, such that 
\meditem{(i)} $\cup G_n$ is everywhere dense in $G$,
\meditem{(ii)} $\mu_n$ are $G_n$-invariant,
\meditem{(iii)} $(X,d,\mu_n)$ form a L\'evy family.
\vskip .3cm
\begfig{3.8cm}
\largeurfig=5cm
\deplacefig=1mm
\montefig=-1mm
\figinsert{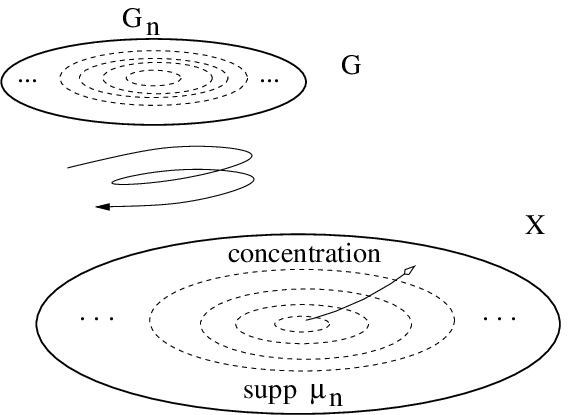}
\figure2{A L\'evy transformation group}
\endfig
In the particular case where $X$ is the group itself equipped with a 
right-invariant metric and the action of $G$ is by left translations,
we say that $G$ is a {\it L\'evy group.}

\parag{Example~10} Let $\H=\ell_2$, and let
$G=\U(\H)_s$, $X=(G,d)$, where $d$ is a right-invariant metric and the 
action is by left translations. Set $G_n=\SU(n)$ (embedded into $\U(\H)$
as a subgroup of block-diagonal operators), 
and let $\mu_n$ denote the normalized Haar
measure on $\SU(n)$. One can view $\mu_n$ as a measure on
all of $\U(\ell_2)_s$ with support $\SU(n)$. 
The $mm$-spaces $(\U(\H)_s,d,\mu_n)$ clearly form a L\'evy family,
because the spaces $(\SU(n)_s,d\vert_{\SU(n)},\mu_n)$ do.
We conclude: $\U(\H)_s$ is a L\'evy group. 

\parag{Example~11}
Let $\pi$ be a strongly continuous
unitary representation of a compact group 
$G$ in $\ell_2$. Then $\ell_2$ decomposes into the
orthogonal direct sum of finite-dimensional (irreducible) unitary $G$-modules,
$\ell_2\cong \bigoplus_{n=1}^\infty V_n$. Set for each $n\in\N$
$$\s_n=\s^{\infty}\cap\bigoplus_{i=1}^n V_n.$$
We obtain a nested sequence of spheres of increasing finite dimension
which are invariant under the action of $G$.
Let $\mu_n$ denote the rotation-invariant probability measure on 
the sphere $\s_n$. Denote also $G_n=G$ for all $n$.
Then $(G,\s^{\infty})$ is a L\'evy
transformation group.

\sec{\\/}{Concentration property and fixed points}
The following definition is an attempt to capture 
`concentration in the absence
of measure' (as indeed there are typically no invariant measures on infinite
dimensional spaces.)

\parag{Definition~7 \cite{M2,M3}}
Let a group $G$ act on a metric space $X$ by uniform isomorphisms. 
Call a subset $A\subseteq X$ {\it essential} if for every 
$\e>0$ and every finite
collection $g_1,\dots,g_N\in G$ one has
$$\bigcap_{i=1}^N g_i A_\e\neq\emptyset.$$ 
(Have another look at Fig. 1!)

\parag{Exercise~5}
The definition obtained by replacing $g_iA_\e$ with $(g_iA)_\e$ is equivalent.
\vskip .3cm
Informally speaking, an essential set is so `big' 
that translates of any $\e$-neighbourhood of it,
taken in any finite number, don't fit in without overlapping.

\parag{Definition~8} ({\it ibid.}) 
A $G$-space $X$ has the {\it concentration property} if
every finite cover of $X$ contains at least one essential set. 
\vskip .3cm
Perhaps one gets a better idea of the property if we start with an
example where it is violated.

\parag{Example~12} (Imre Leader, 1988, unpublished.)
The $\U(\H)$-space $\s^{\infty}$ 
(the unit sphere in $\H=\ell_2$)
does not have the concentration property.
Denote by $E$ the set of all even
natural numbers, and let $P_E$ be
the corresponding projection in $\ell_2$. Set
$$A=\left\{x\in \s^\infty \colon \norm{P_Ex}
\geq {\sqrt 2}/2\right\},$$
$$B=\left\{x\in \s^\infty \colon \norm{P_Ex}
\leq {\sqrt 2}/2\right\}.$$
Clearly, $A\cup B=\s^\infty$. At the same time, both $A$ and $B$ are
inessential. Indeed, let $E_1,E_2,E_3$ be three
arbitrary disjoint infinite
subsets of $\N$, and let $\phi_i\colon \N\to\N$ be bijections
with $\phi_i(E)=E_i$, $i=1,2,3$. Let $g_i$ denote the unitary
operator on $\ell_2(\N)$ induced by $\phi_i$. 
Now
$$g_i(A)=\left\{x\in\s^\infty\colon \norm{P_{E_i}x}\geq
{\sqrt 2} /2\right\},$$
and consequently
$$(g_i(A))_\e\subseteq 
\left\{x\in\s^\infty\colon \norm{P_{E_i}x}\geq
({\sqrt 2} /2)-\e\right\}.$$
Thus, as long as $\e<{\sqrt 2}/2-{\sqrt 3}/3$,
we have
$$\bigcap_{i=1}^3 (g_i(A))_\e=\emptyset.$$
The set $B$ is treated similarly.

\theorem{2}{A compact 
$G$-space $K$ has the concentration property if and only if 
it contains a fixed point: $g\cdot\kappa=\kappa$
for all $g\in G$. 
}

\proof{ $\Rightarrow$: {\it Claim 1.} 
There is a point $\kappa\in K$ such that every neighbourhood of
$\kappa$ is essential.

Assuming the contrary, we could have covered $K$ with inessential open sets
and, selecting a finite open subcover, obtain a contradiction. 

{\it Claim 2.} Any point $\kappa$ as above is $G$-fixed. 

Again, assume that for some $g\in G$, $g\cdot\kappa\neq\kappa$.
Set $\e=d(\kappa,g\cdot\kappa)/2$. 
Choose a number $\delta>0$ so small that $\delta\leq\e/2$ 
and the $g$-translate of
the open ball $B_{\delta}(\kappa)$ is contained in the 
$(\e/2)$-ball around $g\cdot\kappa$.
The set $V=B_\delta(\kappa)$ is essential, yet the 
$\delta$-neighbourhoods of
$V$ and $g\cdot V$ don't meet, a contradiction.

$\Leftarrow$: obvious.
\qed}

The following result provides nontrivial
examples of $G$-spaces with concentration property.

\theorem{3}{Every L\'evy $G$-space $(G,X)$ has the concentration property.
}

\proof{ Let 
$$\gamma=\{A_1,A_2,\dots,A_k\}$$
be a finite cover of $X$.
Since for each $n=1,2,\dots$ the values 
$\mu_n(A_i)$, $i=1,2,\dots,k$, add up
to one, at least one of the sets in $\gamma$, let us denote it simply $A=A_i$,
has the property:
$$\limsup_{n\to\infty} \mu_n(A)\geq\frac 1k.$$
Now let $\e>0$ and a finite collection $g_j,j=1,2,\dots,m$ 
be given.
Using Exercise 2, choose a number $n_0$ so large that 
$$\mu_n (B_\e)>1-\frac 1m$$
whenever $n>n_0$ and $\mu_n(B)\geq\frac 1k$. 
Choose an $n>n_0$ with $\mu_n(A)\geq \frac 1k$; 
then $\mu_n(g_jA)\geq\frac 1k$ as well,
and 
$$\mu_n(g_jA)_\e>1-\frac 1m,~~ i=1,2,\dots,m,$$
implying that the $\e$-neighbourhoods of all the
translates of $A_{\e}$ by $g_j$'s have a common point.
\qed}

To extract useful information from the above, it only remains to link the
concentration property of a $G$-space to that of its compactification.

\lemma{1}{
Let $X$ and $Y$ be two $G$-spaces.\note{As before,
$X$ and $Y$ are metric spaces upon which $G$ acts continuously, by
uniform isomorphisms.}
Let $i\colon X\to Y$ be an equivariant map. If $(G,X)$ has the
concentration property, then so does $(G,Y)$.}

\proof{If $A\subseteq X$ is an
essential subset, then so is $i(A)$. Notice that the uniform continuity of
$i$ is used here in a substantial way.
\qed} 

The following is now immediate.

\theorem{4 \cite{Gr-M1}}{
Let $(G,X)$ be a L\'evy $G$-space and let $K$ be a compact $G$-space,
such that there is an equivariant map $X\to K$. 
Then $K$ has a $G$-fixed point. \qed
}

Using Theorem 1 and Example 10, we obtain

\proclaim{Corollary~1}
{ Whenever the topological group $\U(\ell_2)_s$ continuously acts on a compact
space, it has a fixed point.  
} 

Such topological groups are said to have the {\it fixed point on compacta
property,} or else to be {\it extremely amenable.} 
And indeed, this property is
a drastically strengthened form of the usual amenability, which can be
reformulated as follows (Day): a topological group $G$ 
is amenable if and only if
every affine continuous action of $G$ on a convex compact set [in a locally
convex space] has a fixed point. 

\parag{Remark~5}
No locally compact group can have the fixed point on compacta property,
this is a theorem by Veech (\cite{Ve}, Th. 2.2.1).

\parag{Remark~6} The unitary group $\U(\H)_s$ was the first
`natural' extremely amenable group to be discovered.
The second such
discovery was the group $L_0((0,1),{\Bbb T})$ of all (equivalence
classes of) measurable maps from the unit interval to the circle
rotation group, equipped with the topology of convergence in measure.
This was proved by Glasner (and published years later \cite{Gl})
and, independently, by Furstenberg and Weiss (never published).
This group is a L\'evy group, and the approximating L\'evy family of
subgroups is formed by tori ${\Bbb T}^n$, made up of simple functions
with respect to a refining sequence of measurable
partitions of $(0,1)$. 

It is interesting that both groups mentioned in the previous paragraph
appear as the `outermost' cases of a newly discovered class of
extremely amenable groups. Recall that a von Neumann algebra $M$ is
{\it approximately finite dimensional} if it contains a directed
family of finite-dimensional $\ast$-subalgebras with everywhere dense union.
Denote by $M_\ast$ the predual of $M$. 
It is proved in \cite{G-P} that a von Neumann algebra $M$ is
approximately finite-dimensional if and only if the unitary group
of $M$, equipped with the topology $\sigma(M,M_\ast)$, is
the product of a compact group with an extremely amenable group.

The two cases to consider now are $M={\Cal B}(\H)$,
where the unitary group with the above topology is $\U(\H)_s$, and
$M=L^\infty(0,1)$, in which case the unitary group is 
$L_0((0,1),{\Bbb T})$. 

As a corollary, nuclear $C^\ast$-algebras admit a characterization
in terms of topological dynamics of their unitary groups.
Recall that an action of a group $G$ on a compact space $X$ is
{\it minimal} if the $G$-orbit of every point of $X$ is everywhere dense,
and {\it equicontinuous} if 
the family of all mappings
$x\mapsto gx$, $g\in G$ of $X$ to itself is uniformly equicontinuous.
By considering the enveloping von Neumann algebra, one can deduce that a 
$C^\ast$-algebra $A$ is nuclear if and only if every minimal
continuous action of the
unitary group $\U(A)$, equipped with the $\sigma(A,A^\ast)$-topology,
on a compact space $K$ is equicontinuous.

\parag{Remark~7}
One has to be careful while applying Theorem 4. For instance, 
consider the infinite permutation group $S_\infty$, formed
by all self-bijections of a countably infinite set, say $\Z$.
This group is 
equipped with the natural Polish
topology of pointwise convergence on discrete $\Z$, 
induced by the embedding 
$S_\infty\hookrightarrow\Z^\Z$. 
The idea of applying concentration in finite groups of permutations 
(Example 3) to conclude
that $S_\infty$ is a L\'evy group is attractive, but does not work. 

\parag{Exercise~6} Let $d$ be any right-invariant metric on $S_\infty$,
generating the topology of pointwise convergence. 
Show that $S_\infty$, acting on the left
upon $(S_\infty,d)$, does not have the concentration property. 

[Hint: let $\tau$ be the transposition
exchanging $0$ and $1$ and leaving the rest of $\Z$ fixed. Choose $\e>0$
so that the $\e$-ball around $e_G$ is contained in the intersection
of the isotropy subgroups of $0$ and $1$. Now partition $S_\infty$ into
two sets $A$ and $B$, where
$$
A=\{\sigma\in S_\infty\colon \sigma^{-1}(0)<\sigma^{-1}(1)\}
$$
and $B= S_\infty\setminus A$. Try to apply the 
concentration property to the cover
$\{A,B\}$, the number $\e$, and the collection of
two elements $e,\tau$.]
\vskip .3cm
It follows that $S_\infty$ acts
on some compact space without fixed points.
(This was noted in \cite{P1}.) Very recently such an action was
constructed explicitely by Eli Glasner and Benji Weiss \cite{Gl-W}. 
We will return to their construction later (Subsection 6.4).

\parag{} One can even show that $S_\infty$ is not a L\'evy group
no matter what the group topology is (\cite{P2}, Remark 4.9).
However, it is still possible to put the finite
permutation groups $(S_n)$ 
together so as to obtain a L\'evy group. 

This is the group $\Aut(X,\mu)$ 
of all measure-preserving automorphisms of
the standard non-atomic Lebesgue space, $(X,\mu)$, equipped with 
the weak topology, that is,
the weakest topology making every map of the form
$\Aut(X,\mu)\ni \tau\mapsto \mu(A\,\Delta\,\tau(A))\in\R$ continuous, where
$A\subseteq X$ is a measurable set.
This group contains finite permutation groups, realized as subgroups of
interval exchange transformations, and any right-invariant metric
makes those subgroups into a L\'evy family. 
A similar result holds for the group $\Aut^\ast(X,\mu)$ of all measure class
preserving transformations.
(Thierry Giordano and the author, \cite{G-P}). 

\sec{\\/}{Invariant means on spheres}

Let a group $G$ act on a metric space $X$
by uniform isomorphisms. The formula
$$^gf(x)=f(g^{-1}\cdot x)$$
determines an action of $G$ on the space $\UCB(X)$ of all uniformly continuous
bounded complex valued functions on $X$ by linear isometries.
If $G$ is a topological group acting on $X$ continuously,
the above action of $G$ on $\UCB(X)$ need not, in general, be continuous.
(An example: $G=\U(\ell_2)_s$, $X=\s^\infty$.)
However, the action will be continuous if $X$ is compact. (An easy check.)
To some extent, the latter observation can be inverted.

\parag{Exercise~7}
Let a topological group $G$ act continuously on a commutative 
unital $C^\ast$-algebra
$A$ by automorphisms. Then this action determines a continuous
action of $G$ on the space of maximal ideals of $A$, 
equipped with the usual (weak$^\ast$) topology.

\parag{}
Recall that a {\it mean} on a space $\Cal F$ of functions is a
positive linear functional, $m$, of norm one, sending the function $1$ to $1$.
A mean is {\it multiplicative} if $\Cal F$ is an algebra
and the mean is a homomorphism of this algebra to $\C$.

\proclaim{Corollary~2}
{ Let $(G,X)$ be a L\'evy $G$-space. Then there exists a $G$-invariant
multiplicative
mean on the space $\UCB(X)$ of all bounded uniformly 
continuous functions on $X$. 
} 

\proof{ According to Exercise 7, the group
$G$ acts continuously on the space $\frak M$ of maximal ideals of the
$C^\ast$-algebra $\UCB(X)$.
Therefore, $\frak M$ is an equivariant
compactification of $X$.
By force of Theorem 4, there is a fixed point $\phi\in{\frak M}$, 
which is the desired invariant multiplicative mean.
\qed}

The following is deduced by considering Example 11.

\proclaim{Corollary~3~{\rm\cite{Gr-M1}}}
{ If a compact group $G$ is represented by unitary operators in
an infinite-dimensional Hilbert space
$\H$, then there exists a $G$-invariant multiplicative mean on the
uniformly continuous bounded functions on the unit sphere of $\H$.
}

\parag{Remark~8} The infinite-dimensionality of $\H$ is essential. 
Since the unit sphere $\s$ of a finite-dimensional space $\H$ is compact,
an invariant multiplicative mean on $\UCB(\s)$ 
exists if and only if there is a fixed vector $\xi\in\s$.
\vskip .3cm
Means on $\UCB(X)$, where $X=\s^{\infty}$ is the unit sphere in the Hilbert
space, as well as 
some other infinite dimensional manifolds, were studied by
Paul L\'evy, who viewed them as (substitutes for)
infinite-dimensional integrals.\note{The multiplicativity of
some of those means, which is not exactly a property 
one expects of an integral,
becomes clear if one recalls an equivalent way to express the concentration
phenomenon: on a high-dimensional structure, every 1-Lipschitz function is,
probabilistically, almost constant, cf. Section 7.}
 The invariant means
can thus serve as a substitute for invariant integration 
on the infinite-dimensional
spheres. One can substantially generalize Corollary 3.
With this purpose in view, it is convenient to enlarge
the concept of a L\'evy transformation group.

If $\mu_1,\mu_2$ are probability measures on the same metric space
$X$, then the
{\it transportation distance} between them is defined as
$$d_{tran}(\mu_1,\mu_2) =
\inf \int_{X\times X} d(x,y\,) d\nu(x,y),$$
where the infimum is taken over all probability measures $\nu$ on the
product space $X\times X$ such that $(\pi_i)_\ast\nu=\mu_i$
for $i=1,2$ and $\pi_1,\pi_2\colon X\times X\to X$ denote the
coordinate projections.

The way to think of the transportation distance is to identify 
each probability measure with a pile of sand, then 
$d_{tran}(\mu_1,\mu_2)$ is the minimal average
distance that each grain of sand has to travel when the first pile
is being moved to take place of the second.\note{In computer
science, the transportation distance is known as the
Earth Mover's Distance (EMD).}

Let us from now on replace Definition 6 with the following,
more general one.

\parag{Definition~9}
Say that a $G$-space $(G,X)$ is L\'evy if there is a 
net of probability measures
$(\mu_\alpha)$ on $X$, such that the $mm$-spaces
$(X,d,\mu_\alpha)$ form a L\'evy family and for each $g\in G$ 
$$
d_{tran}(\mu_\alpha,g\mu_\alpha)\to 0.
$$

Theorems 3 and 4 
remain true, with very minor modifications of the
proofs. 

Here is one application.
A unitary representation $\pi$ of a group $G$ in a
Hilbert space $\H$ is {\it amenable} in the sense of Bekka 
\cite{Be} if there exists a 
state, $\phi$, on the algebra ${\Cal B}(\H)$
of all bounded operators on the space $\H$ of representation, which is
invariant under the action of $G$ by inner automorphisms: 
$\phi(\pi_gT\pi_g^\ast)=\phi(T)$ for
every $T\in B(\H)$ and every $g\in G$. 

\theorem{5 \cite{P2}}{
Let $\pi$ be a unitary representation of a group $G$ in a Hilbert space $\H$.
The following are equivalent.
\meditem{\rm (i)} $\pi$ is amenable.
\meditem{\rm (ii)} Either $\pi$ has a finite-dimensional subrepresentation, or 
$(G,\s)$ has the concentration property {\rm(}or both{\rm)}. 
\bigitem{\rm (iii)} There is a $G$-invariant mean on the space $\UCB(\s)$
{\rm(}a `L\'evy-type integral.'\,{\rm)}
} 

\proof{ (i) $\Rightarrow$ (ii): according to
Th. 6.2 and Remark 1.2.(iv) in \cite{Be}, a representation $\pi$
is amenable if and only if for every finite set 
$g_1,g_2,\dots,g_k$ of elements of
$G$ and every $\e>0$ there is a projection $P$ of finite rank such that for all
$i=1,2,\dots,k$
$$\norm{P - \pi_{g_i} P \pi_{g_i}^\ast}_1<\e \norm P_1,$$
where $\norm\cdot_1$ denotes the trace class operator norm. 
It follows that the transportation distance between the
Haar measure on the unit sphere in the range of the
projection $P$ and the translates of this measure by operators
$\pi_{g_i}$ can be made as small as desired
via a suitable choice of $P$. 
Now a variant of Theorem 4 applies. (See \cite{P2} for details.) 
\par
(ii) $\Rightarrow$ (iii): in the first case, the mean is obtained by
invariant integration on the finite-dimensional 
sphere, while in the second case even a multiplicative mean
exists. \par
(iii) $\Rightarrow$ (i): 
Let $\psi$ be a $G$-invariant mean on $\UCB(\s_\H)$.
For every bounded linear operator $T$ on
$\H$ define a (Lipschitz) function $f_T\colon\s_\H\to\C$ by
$$
\s_\H\ni\xi\mapsto f_T(\xi):=\langle T\xi, \xi\rangle\in \C,
$$
and set $\phi(T):=\psi(f_T)$. This $\phi$ is a $G$-invariant
mean on ${\Cal B}(\H)$.
\qed}

\proclaim{Corollary~4}{
A locally compact group $G$ is amenable if and only if for every strongly
continuous unitary representation of $G$ in an infinite-dimensional Hilbert
space the pair $(G,\s^\infty)$ has the property of concentration.
}

\proclaim{Corollary~5}{ There is no invariant mean 
on $\UCB(\s^{\infty})$ for the
full unitary group $\U(\ell_2)$.
}

\proof{If such a mean existed, then every unitary representation of
every group
would be amenable, in particular every group would be amenable
(by Th. 2.2 in \cite{Be}).}

(Of course Corollary 5 also follows from Imre Leader's Example 12 modulo
Theorem 2 and Lemma 1.)

A (not necessarily locally compact) topological group $G$ is 
{\it amenable} if there is a left-invariant mean on the space
${\operatorname{RUCB}}(G)$ of all right uniformly continuous bounded
functions on $G$. 
Denote by $\U(\ell_2)_u$ the full unitary group with the 
uniform operator topology. 

\proclaim{Corollary~6}{{\rm [Pierre de la Harpe \cite{dlH},
proved by different means]}
The topological group $\U(\ell_2)_u$ is not amenable. 
}

\proof{Choose an arbitrary $\xi\in\s^\infty$. To every
function $\psi\in\UCB(\s^\infty)$ associate the function
$\tilde\psi$ as follows:
$$
G\ni g\mapsto \tilde \psi(g) := \psi(\pi_g^\ast(\xi))\in \C.
$$
The correspondence $\psi\mapsto\tilde\psi$ is a $G$-equivariant 
positive bounded unit-preserving linear operator
from $\UCB(\s^\infty)$ to ${\operatorname{RUCB}}(\U(\ell_2)_u)$,
and any left-invariant mean $\varphi$ on the latter $G$-module 
would thus determine a $G$-invariant mean on the former
$G$-module, contradicting Corollary 5.
\qed}

\parag{Example~13}
In a similar fashion, by considering the action of $\Aut(X,\mu)$ on
$L^2_0(X,\mu)$, where 
$X={\operatorname{SL}}(3,\R)/{\operatorname{SL}}(3,\Z)$, one
deduces that $\Aut(X,\mu)_u$ with the uniform topology is not amenable
\cite{G-P}.

\sec{\\/}{Ramsey--Dvoretzky--Milman property}

\subsec{\\/}{Extreme amenability and small oscillations}
One way to intuitively describe a `Ramsey-type result' is as follows.
Suppose $\frak X$ is a large (and often highly homogeneous)
structure of some sort or other. 
Let $\frak X$ be partitioned into a finite number of
pieces in an arbitrary way. No matter how irregular and `ragged' the
pieces are, at least one of them always contains the remnants of
the original structure, that is, a (possibly much smaller, but still 
detectable) substructure of the same type which survived
intact.

We are now going to explicitely link the fixed point on compacta
property to Ramsey-type results. 
Here is the first step.

\parag{Exercise~8} Prove that a topological group $G$ is extremely amenable
if and only if for every finite collection $g_1,\dots,g_n$ of
elements of $G$, every bounded right uniformly continuous function
$f\colon G\to \R^N$ from $G$ to a finite-dimensional Euclidean space,
and every $\e>0$ there is an $h\in G$ such that
$\abs{f(h)-f(g_ih)}<\e$ for each $i=1,2,\dots,n$. 
\vskip .3cm

[{\it Hints:} $\Rightarrow$: the action of $G$ on the space 
${\Cal S}(G)$ of
maximal ideals of the $C^\ast$-algebra ${\operatorname{RUCB}}(G)$
is continuous, and $G$ itself can be thought of as an everywhere
dense subset of ${\Cal S}(G)$. \vskip .2mm
$\Leftarrow$: 
form a net of suitably indexed elements $h$ as above and
consider any limit point of the net $h_\alpha\cdot\xi$, where $\xi$
is an arbitrary element of the compact space upon which $G$ acts
continuously.]

\parag{Exercise~9} Prove that the above condition for extreme
amenability is, in turn, equivalent to the following.
For every bounded
{\it left} uniformly continuous function $f$ from $G$ to
a finite-dimensional Euclidean space, every finite subset $F$ of $G$,
and every $\e>0$, the oscillation of
$f$ on a suitable left translate of $F$ is less than $\e$:
$$
\exists g\in G,~~\operatorname{Osc}(f\vert_{gF})<\e.
$$

It is convenient to deal with the above property in a more
general context of $G$-spaces.

\parag{Definition~11 {\rm \cite{Gr1}}}
Say that a $G$-space $X$ (in our agreed sense)
has the {\it Ramsey--Dvoretzky--Milman property} if for every
bounded
uniformly continuous function $f$ from $X$ to a finite-dimensional
Euclidean space, every $\e>0$, and every finite $F\subseteq X$,
there is a $g\in G$ with the property
$$
\operatorname{Osc}(f\vert_{gF})<\e.
$$

\begfig{3.4cm}
\largeurfig=5cm
\deplacefig=1mm
\montefig=-1mm
\figinsert{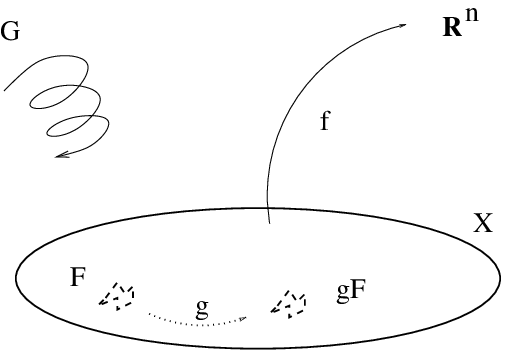}
\figure3{The Ramsey--Dvoretzky--Milman property}
\endfig
\parag{Remark~9} Equivalently, $F$ can be assumed compact.

\proclaim{Corollary~7}{For a topological group $G$ the following
are equivalent:
\item{\rm (i)} $G$ is extremely amenable,
\item{\rm (ii)} every metric space $X$ upon which $G$ acts continuously
and transitively by isometries has the R--D--M property,
\item{\rm (iii)} every homogeneous factor-space $G/H$, equipped with a 
left-invariant metric {\rm(}or the left uniform structure{\rm)}, has  
the R--D--M property.}

Next, we will discover two very important situations where the
R--D--M property appears naturally.

\subsec{\\/}{Dvoretzky theorem.}
Here is the famous result.

\theorem{(Arieh Dvoretzky)}{For all $\e>0$ there is a constant
$c=c(\e)>0$ such that for any $n$-dimensional normed space
$(X,\norm\cdot_E)$ there is a subspace $V$ of $\dim V\geq c\,\log n$ and a
Euclidean norm $\norm\cdot_2$
with $\norm x_2\leq\norm x_E\leq (1+\e)\norm x_2$ for all $x\in V$.
}

The studies of the phenomenon of concentration of measure were 
given a boost by Vitali Milman's new proof of the 
Dvoretzky theorem \cite{M1},
based on a suitable finite-dimensional approximation to the lemma which
directly follows from results that we have previously stated:

\lemma{(Milman)}{The pair $(\U(\H),\s^{\infty})$ has the R--D--M property,
where $\s^{\infty}$ is the unit sphere of an infinite-dimensional Hilbert space
$\H$.}

\subsec{\\/}{Ramsey theorem}
Let $r$ be a positive natural number. By $[r]$ one denotes the
set $\{1,2,\dots,r\}$. 
A {\it colouring} of a set $X$ with $r$ colours, or simply
$r${\it -colouring,} is any map
$\chi\colon X\to [r]$. 
A subset $A\subseteq X$ is {\it monochromatic}
if for every $a,b\in A$ one has $\chi(a)=\chi(b)$. 

Put otherwise, a finite colouring of a set $X$ is nothing but a partition
of $X$ into finitely many (disjoint) subsets. 

Let $X$ be a set, and let $k$ be a natural number. Denote by
$[X]^k$ the set of all $k$-subsets of $X$, that is,
all (unordered!) subsets containing exactly $k$ elements. 

\proclaim{Infinite Ramsey Theorem}{
Let $X$ be an infinite set, and let $k$ be a natural number.
For every finite colouring of $[X]^k$ there exists an infinite
subset $A\subseteq X$ such that the set $[A]^k$ is monochromatic.}

\begfig{4.3cm}
\largeurfig=7cm
\deplacefig=7mm
\montefig=-1mm
\figinsert{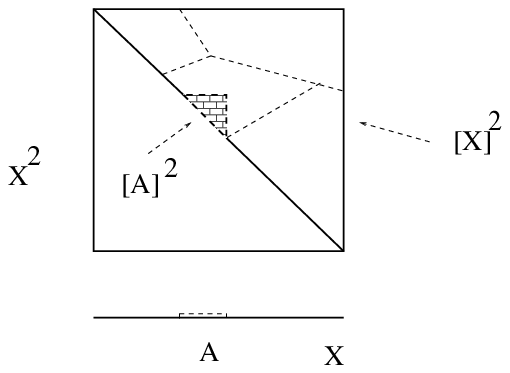}
\figure4{Ramsey theorem for $k=2$}
\endfig

\parag{Remark~10}
For $k=1$ the statement is simply the pigeonhole principle.
Here is a popular interpretation of the result in the case $k=2$.
Among infinitely many people, either there is an infinite subset
of people every two of whom know each other, or there is an
infinite subset no two members of which know each other.

\proclaim{Finite Ramsey theorem}{
For every triple of natural numbers, $k,l,r$, there exists a natural
number $R(k,l,r)$ with the following property.
If $N\geq R(k,l,r)$ and the set of all $k$-subsets of
$[N]$ is coloured using $r$ colours, then there is a subset
$A\subseteq [N]$ of cardinality $\abs A= l$ such that all $k$-subsets
of $A$ have the same colour.}

\parag{Remark~11}
The Infinite Ramsey Theorem implies the finite version through a
simple compactness argument.
At the same time, the infinite version does not seem to quite
follow from the finite one. The finite version is equivalent
to the following statement:
\vskip .3cm
{\it Let $X$ be an infinite set, and let $k$ be a natural number.
For every finite colouring of $[X]^k$ and every natural $n$ 
there exists a
subset $A\subseteq X$ of cardinality $n$
such that $[A]^k$ is monochromatic.}
\vskip .3cm
A good introductory reference to Ramsey theory is \cite{Gra}.
\vskip .3cm
Denote by $\operatorname{Aut}(\Q)$ the group of all
order-preserving bijections of the set of rational numbers, equipped
with the topology of pointwise convergence on the discrete set $\Q$.
In other words, we regard $\operatorname{Aut}(\Q)$ as a (closed)
topological subgroup of $S_\infty$. 
A basic system of neighbourhoods of identity is formed by open subgroups
each of which
stabilizes elements of a given finite subset of $\Q$. 

\parag{Exercise~10} Use Corollary 7 to 
prove that the finite Ramsey theorem is 
equivalent to the statement: 
\vskip .2cm
{\it the topological group $\operatorname{Aut}(\Q)$ is extremely
amenable.}
\vskip .3cm
[{\it Hint:} for a finite subset $M\subset\Q$, the left factor space
of $\operatorname{Aut}(\Q)$ by the stabilizer of $M$ can be identified
with the set $[\Q]^{n}$, where $n=\abs M$,
equipped with the discrete uniformity (or $\{0,1\}$-valued metric).
Cover $[\Q]^{n}$ with finitely many sets on each of which the given
function $f$ has oscillation $<\e$, and apply Ramsey theorem.
Use Remark 11.]

\subsec{\\/}{Extreme amenability and minimal flows}

\proclaim{Corollary~8}{The group of orientation-preserving
homeomorphisms of the closed
unit interval, $\operatorname{Homeo}_+(\I)$, equipped with the compact-open
topology, is extremely amenable.}
\proof{
Indeed, the extremely amenable
group $\operatorname{Aut}(\Q)$ admits a continuous
monomorphism with a dense image into the group
$\operatorname{Homeo}_+(\I)$.}

\parag{Remark~12} Thompson's group $F$ consists of all piecewise-linear
homeomorphisms of the interval whose points of
non-smoothness are finitely many dyadic rational numbers, 
and the slopes of any linear part
are powers of 2. (See \cite{CFP}.)
It is a major open question in combinatorial group theory whether the
Thompson group is amenable. Since $F$ is everywhere dense in 
$\operatorname{Homeo}_+(\I)$, 
our Corollary 8 does not contradict the possible
amenability of $F$.
\vskip .3cm

Using
extreme amenability of the topological groups $\Aut(\Q)$ and
$\operatorname{Homeo}_+(\I)$, one is able
to explicitely compute the
universal minimal flows of some larger topological groups as follows.

\proclaim{Corollary~9}
{The circle $\s^1$ forms the universal minimal
\linebreak
$\operatorname{Homeo}_+(\s^1)$-space.}

\proof{Let $\theta\in\s^1$ be an arbitrary point.
The isotropy subgroup 
$\operatorname{St}_{\theta}$ of $\theta$ is isomorphic to
$\operatorname{Homeo}_+(\I)$. Because of that,
whenever the topological group $\operatorname{Homeo}_+(\s^1)$
acts continuously on a compact space $X$, the subgroup
$\operatorname{St}_{\theta}$ has a fixed point,
say $x'\in X$. The mapping  $\operatorname{Homeo}_+(\s^1)\ni h\mapsto h(x')
\in X$ is constant on the left $\operatorname{St}_{\theta}$-cosets and
therefore gives rise to a continuous equivariant map 
$\operatorname{Homeo}_+(\s^1)/\operatorname{St}_{\theta}\cong\s^1\to X$. 
}

For the above results concerning groups
$\operatorname{Aut}(\Q)$,
$\operatorname{Homeo}_+(\I)$, and $\operatorname{Homeo}_+(\s^1)$,
see \cite{P1}.
\vskip .3cm
Now denote by $\LO$ the set of all linear orders on $\Z$,
equipped with the (compact) 
topology induced from $\{0,1\}^{\Z\times\Z}$. The group $S_\infty$
acts on $\LO$ by double permutations.

\parag{Exercise~11} Prove that the action of $S_\infty$ on
$\LO$ is continuous and minimal
(that is, the orbit of each linear order is everywhere dense
in $\LO$). 

\parag{} Recall that a linear order $\prec$ is called {\it dense}
if it has no gaps. A dense linear order without  
least and greatest elements is said to be of type $\eta$.
The collection $\LO_\eta$ of all
linear orders of type $\eta$ on $\Z$
can be identified with the 
factor space $S_\infty/\Aut(\prec)$ through the correspondence
$\sigma\mapsto\,^\sigma\!\!\prec$. Here $\prec$ is some chosen
linear order of type $\eta$ on $\Z$ and $\Aut(\prec)$ stands for
the group of order-preserving self-bijections of $(\Z,\prec)$,
acting on the space of orders in a natural way:
$(x\,\,^\sigma\!\!\prec y)\Leftrightarrow \sigma^{-1}x\prec \sigma^{-1}y$.

\parag{Exercise~12} Show that under the above identification 
the uniform
structure on $\LO_\eta$, induced from the compact space $\LO$, is
the finest uniform structure making the quotient map
$S_\infty\to S_\infty/\Aut(\prec)\cong \LO_\eta$ right uniformly
continuous.

\parag{} Let now $X$ be a compact $S_\infty$-space.
The topological subgroup $\Aut(\prec)$ of $S_\infty$ has a fixed point
in $X$, say $x'$. (Exercise 10.) The mapping  
$S_\infty\ni\sigma\mapsto  \sigma(x')\in X$ is constant on the 
left $\Aut(\prec)$-cosets and thus gives rise to 
a mapping $\phi\colon\LO_\eta\to X$.
Using Exercise 12, it is easy to see that 
$\phi$ is right uniformly continuous and thus extends to a morphism
of $S_\infty$-spaces $\LO\to X$. We have established the following
result.

\theorem{6 (Glasner and Weiss \cite{Gl-W})}
{The compact space $\LO$ forms the universal minimal $S_\infty$-space.}

\subsec{\\/}{The Urysohn metric space}

The {\it universal Urysohn metric space} $\Ur$ \cite{Ur}
is determined uniquely
(up to an isometry) by the following conditions:
\vskip .2cm
\item{(i)} $\Ur$ is a complete separable metric space;
\item{(ii)} $\Ur$ is $\omega$-homogeneous, that is, every isometry
between two finite subspaces of $\Ur$ extends to an isometry of $\Ur$;
\item{(iii)} $\Ur$ contains an isometric copy of every separable metric
space.
\vskip .3cm
A probabilistic description of this space was given by Vershik
\cite{Ver}:
the completion of the space of integers equipped with a
`sufficiently random' metric is almost surely isometric to $\Ur$.

The group of isometries $\Iso(\Ur)$ with the compact-open topology
is a Polish (complete metric separable) topological group, 
which also possesses a universality property: it contains an
isomorphic copy of every separable metric group \cite{Usp}.
See also \cite{Gr3}.

Using concentration of measure, one can prove that the 
group  $\Iso(\Ur)$ is extremely amenable. The Ramsey--Dvoretzky--Milman
property leads to the following Ramsey-type result:
\vskip .3cm
{\it 
Let $F$ be a finite metric space, and
let all isometric embeddings of $F$ into $\Ur$
be coloured using finitely many colours. 
Then for every finite metric space $G$ and every $\e>0$ 
there is an isometric copy $G'\subset \Ur$ of $G$ such that
all isometric embeddings of $F$ into $\Ur$ that factor through
$G$ are monochromatic to within $\e$.}

\begfig{3.6cm}
\largeurfig=5cm
\deplacefig=1mm
\montefig=-1mm
\figinsert{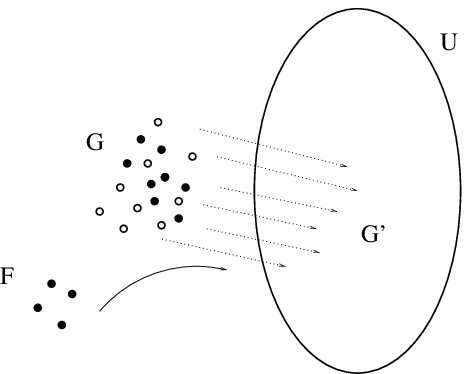}
\figure5{A Ramsey-type result for metric spaces}
\endfig

Here we say that a set $A$ is {\it monochromatic to within $\e$}
if there is a monochromatic set $A'$ at a Hausdorff distance $<\e$
from $A$. In our case, the Hausdorff distance is formed with 
regard to the uniform metric on $\Ur^F$.

One can also obtain similar results, for example, for the separable
Hilbert space $\ell_2$ and for the unit sphere $\s^{\infty}$ in $\ell_2$
\cite{P3}.

\sec{\\/}{Concentration to a non-trivial space} 

Let $f$ be a Borel measurable real-valued function on an $mm$-space
$X=(X,d,\mu)$. 
A number $M=M_f$ is called a {\it median} (or {\it L\'evy mean})
of $f$ if both $f^{-1}[M,+\infty)$ and
$f^{-1}(-\infty,M]$ have measure $\geq\frac 12$. 

\parag{Exercise~13} Show that the median $M_f$ always exists, though need
not be unique.

\parag{Exercise~14} Assume that a function $f$ as above is 1-Lipschitz,
that is, $\vert{f(x)-f(y)}\vert\leq d(x,y)$ for all
$x,y\in X$. Prove that for every $\e>0$
$$
\mu\{\vert f(x)-M_f\vert>\e\}\leq 2\alpha_X(\e).
$$

\parag{}
Thus, one can express the phenomenon of concentration of measure by
stating that on a `high-dimensional' $mm$-space, every Lipschitz
(more generally, uniformly continuous) function is, probabilistically,
almost constant.

Following Gromov \cite{Gr3, 3$\frac 12$.45}, let us recast
the concentration phenomenon yet again.

On the space $L(0,1)$ of all measurable functions define the metric
$\operatorname{me}_1$, generating the topology of convergence in measure,
 by letting 
$\operatorname{me}_1(h_1,h_2)$ stand for the infimum of all $\lambda>0$ with
the property 
$$
\mu^{(1)}\{\vert h_1(x)-h_2(x)\vert>\lambda\}<\lambda.
$$
(Here $\mu^{(1)}$ denotes the
Lebesgue measure on the unit interval $\I=[0,1]$.)

Now let $X=(X,d_X,\mu_X)$ and
$Y=(Y,d_Y,\mu_Y)$ be two Polish $mm$-spaces. 
There exist measurable maps
$f\colon\I\to X$, $g\colon\I\to Y$ such that $\mu_X=f_\ast\,\mu^{(1)}$
and $\mu_Y=g_\ast\,\mu^{(1)}$. Denote by $L_f$
the set of all functions of the form $h= h_1\circ f$, where
$h_1\colon X\to \R$ is 1-Lipschitz, having the property $h(0)=0$.
Similarly, define the set $L_g$. 
Now define a non-negative real number
$\HL(X,Y)$ as the infimum of Hausdorff distances between $L_f$
and $L_g$ (formed using the metric $\operatorname{me}_1$ on the space of
functions), taken over all parametrizations $f$ and $g$ as above.

\begfig{4.4cm}
\largeurfig=5cm
\deplacefig=1mm
\montefig=-1mm
\figinsert{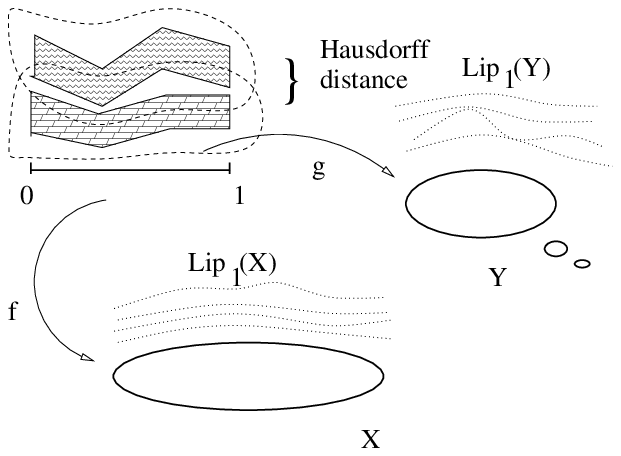}
\figure6{Gromov's distance $\HL$ between $mm$-spaces}
\endfig

\parag{Exercise~15}
Prove that $\HL$ is a metric on the space of (isomorphism classes
of) all Polish $mm$-spaces.

\parag{Exercise~16} Prove that a sequence of $mm$-spaces
$X_n=(X_n,d_n,\mu_n)$ forms a L\'evy family if and only if
it converges to the trivial $mm$-space in the metric $\HL$:
$$
X_n\overset{\HL}\to{\longrightarrow} \{\ast\}.
$$

If one now replaces the trivial space on the right hand side
with an arbitrary $mm$-space,\note{Or, more generally,
a uniform space --- for instance, a non-metrizable compact
space --- with measure.}
one obtains the concept of
{\it concentration to a non-trivial space.}

According to Gromov, this type of concentration commonly occurs
in statistical physics. At the same time, there are very few known
non-trivial examples
of this kind in the context of transformation
groups.

Here is just one problem in this direction, suggested by Gromov.
Every probability measure $\nu$ on a group $G$ determines
a random walk on $G$. 
How to associate to $(G,\nu)$ in a natural way 
a sequence of $mm$-spaces which
would concentrate to the boundary \cite{Fur}
of the random walk? 

\sec{\\/}
{Reading suggestions} The 2001 Borel seminar was based on the
Chapter $3\frac 12$ of the green book \cite{Gr3}, 
which contains a wealth of ideas and concepts and
may be complemented by \cite{Gr4}.
The survey \cite{M3} by Vitali Milman, to whom we owe the present status of
the concentration of measure phenomenon,
is highly relevant and rich in material, 
especially if read in
conjunction with a recent account of the subject by the same author
\cite{M4}.
The book \cite{M-S} is, in a sense, indispensable and 
should always be within one's reach. 
Talagrand's fundamental paper \cite{Ta1} has to be at least browsed by every
learner of the subject, while the paper \cite{Ta2} of the same
author offers an
independent introduction in the subject of concentration of measure.
The newly-published book by Ledoux \cite{Led}, apparently the
first ever monograph devoted exclusively to concentration, is highly readable
and covers a wide range of topics. 
Don't miss the introductory survey by Schechtman
\cite{Sch}.
The modern setting for concentration was designed
in the important paper \cite{Gr-M1} by Gromov and Milman, 
which had also introduced the
subject of this lecture and from where many results (perhaps with slight
modifications) have been taken.

\begreflab{References}{1.4cm}
 
\ref A-M
{\smc Alon, N. {\rm and} V.D. Milman.}
$\lambda_1$, isoperimetric inequalities for graphs, and
superconcentrators.
{\it J. Comb. Theory Ser. B 38\/} (1985), 73--88.

\ref Be
{\smc Bekka,  M.E.B.}
Amenable unitary representations of locally compact groups.
{\it Invent. Math. 100\/} (1990), 383--401.

\ref{CFP}
{\smc Cannon, J.W., W.J. Floyd, {\rm and} W.R. Parry.}
Introductory notes on Richard Thompson's groups. 
{\it Enseign. Math. (2) 42\/} (1996), 215--256. 

\ref{Fur}
{\smc Furstenberg, H.}
Random walks and discrete subgroups of Lie groups, in:
{\it Advances in Probability and Related Topics} 
(P. Ney, ed.), Vol. I, pp. 1--63, Marcel Dekker,
New York, 1971.


\ref{G-P} {\smc Giordano, T. {\rm and} V. Pestov.}
Some extremely amenable groups. 
{\it C.R. Acad. Sc. Paris, S\'er. I, 334\/} (2002), No. 4,
273--278.

\ref{Gl} {\smc Glasner, S.}
On minimal actions of Polish groups. {\it Topology Appl. 85\/}
(1998), 119--125.

\ref{Gl-W} {\smc Glasner, E. {\rm and} B. Weiss.}
Minimal actions of the group ${\Bbb S}({\Bbb Z})$ of
permutations of the integers.
{\it Geom. Funct. Anal. 12\/}  (2002), 964--988. 

\ref{Gra} {\smc Graham, R.L.}
{\it Rudiments of Ramsey theory.}
Regional Conference Series in Mathematics 45,
American Mathematical Society, Providence, R.I., 1981.

\ref{Gr1} {\smc Gromov, M.}  Filling Riemannian manifolds.
{\it J. Differential Geom. 18\/} (1983), 1--147.

\ref{Gr2} {\smc Gromov, M.} 
Isoperimetric inequalities in Riemannian manifolds,
Appendix I in \cite{M-S}, 114--129.

\ref{Gr3} {\smc Gromov, M.} 
{\it Metric Structures for Riemannian and
Non-Riemannian Spaces.} Birkh\"auser Verlag, 1999.

\ref{Gr4} {\smc Gromov, M.} Spaces and questions. GAFA 2000 (Tel Aviv, 1999).
{\it Geom. Funct. Anal. Special Volume}, 
 Part I (2000), 118--161.

\ref{Gr-M1}
{\smc Gromov, M. {\rm and} V.D. Milman.}
A topological application of the isoperimetric inequality.
{\it Amer. J. Math. 105\/} (1983), 843--854.

\ref{Gr-M2}
{\smc Gromov, M. {\rm and} V.D. Milman.}
Generalization of the spherical isoperimetric inequality to
uniformly convex Banach spaces.
{\it Compositio Math. 62\/} (1987), 263--282.

\ref{dlH} {\smc de la Harpe, P.}
Moyennabilit\'e de quelques groupes
topologiques de dimension infinie.
{\it C.R. Acad. Sc. Paris, S\'er. A 277\/} (1973), 1037--1040.


\ref{Led} {\smc Ledoux, M.}
{\it The Concentration of Measure Phenomenon,}
Mathematical Surveys and Monographs 89, American Mathematical
Society (Providence), 2001.

\ref{L\'ev} 
{\smc L\'evy, P.} {\it Le\c cons d'analyse fonctionnelle,} 
Gauthier-Villars (Paris), 1922.

\ref{Ma} {\smc Maurey, B.} 
Constructions de suites sym\'etriques.
{\it C.R. Acad. Sci. Paris, S\'er. A-B 288\/} (1979), 679--681.


\ref{M1} {\smc Milman, V.D.} A new proof of the theorem of
A. Dvoretzky on sections of convex bodies.
{\it Functional Anal. Appl. 5\/} (1971), no. 4, 288--295

\ref{M2} {\smc Milman, V.D.} Diameter of a minimal invariant subset
of equivariant Lipschitz actions on compact subsets of
$\R^k$. In: {\it Geometrical Aspects of Functional Analysis
(Israel Seminar, 1985--86).\/} Lecture Notes in Math.  1267
(1987), Springer-Verlag (Berlin), pp. 13--20.

\ref{M3}  {\smc Milman, V.D.} 
The heritage of P. L\'evy in
geometrical functional analysis. {\it Ast\'erisque 157--158\/}
(1988), 273--301.

\ref{M4} {\smc Milman, V.D.} 
Topics in asymptotic geometric analysis. GAFA 2000 (Tel Aviv, 1999).
{\it Geom. Funct. Anal. Special Volume},
 Part I (2000), 792--815.

\ref{M-S} {\smc Milman, V.D. {\rm and} G. Schechtman.}
{\it Asymptotic Theory of Finite Dimensional Normed Spaces.}
Lecture Notes in Math. 1200, Springer-Verlag, 1986.


\ref{P1} {\smc Pestov, V.G.}
On free actions, minimal flows, and a problem by Ellis.
{\it Trans. Amer. Math. Soc. 350\/} (1998), 
pp. 4149--4165.

\ref{P2} {\smc Pestov, V.G.}
Amenable representations and dynamics of the unit sphere
in an infinite-dimensional Hilbert space.
{\it Geom. Funct. Anal. 10\/} (2000), 1171--1201.

\ref{P3} {\smc Pestov, V.} 
Ramsey--Milman phenomenon, Urysohn metric spaces,
and extremely amenable groups. {\it Israel J. Math. 127\/}
(2002), 317--358.

\ref{Sch} {\smc Schechtman, G.} 
Concentration, results and applications. 
\newline E-print
{\tt http://www.wisdom.weizmann.ac.il/home/gideon/}\newline
{\tt public$\_$html/recentPubs.html}



\ref{Ta1} {\smc Talagrand, M.} 
Concentration of measure and
isoperimetric inequalities in product spaces. 
{\it Inst. Hautes \'Etudes Sci. Publ. Math.
81\/} (1995), 73--205.

\ref{Ta2}
{\smc Talagrand, M.}
A new look at independence. {\it Ann. Probab. 24\/} (1996), 1--34. 

\ref{Ta3} {\smc Talagrand, M.} 
New concentration inequalities in product spaces.
{\it Invent. Math. 126\/} (1996), 505--563. 

\ref{Ur}
{\smc Urysohn, P.} Sur un espace m\'etrique universel.
{\it Bull. Sci. Math. 51\/} (1927), 43--64, 74--90.

\ref{Usp}
{\smc Uspenskij, V.V.}  
On the group of isometries of the
Urysohn universal metric space.
{\it Comment. Math. Univ. Carolinae 31\/}
(1990), 181--182.

\ref{Ve} {\smc Veech, W.A.} Topological dynamics.
{\it Bull. Amer. Math. Soc. 83\/} (1977), 775--830.

\ref{Ver}
{\smc Vershik, A.M.} The universal Urysohn space, 
Gromov metric triples and
random metrics on the natural numbers.
{\it Russian Math. Surveys 53\/} (1998), 921--928;
corrigendum, {\it ibid. 56} (2001), p. 1015.

\endref
 
\adresse
Vladimir Pestov \\
Department of Mathematics and Statistics \\
University of Ottawa \\
585 King Edward Ave. \\
Ottawa, ON K1N 6N5 \\
Canada
{\it e-mail:\/} vpest283$\@$science.uottawa.ca
\endadresse
 
\bye